\documentclass[12pt]{article}

\usepackage{mathptmx,amsmath,amssymb,amsfonts}
\usepackage{amsthm} 
\usepackage{graphicx,color,xcolor,pictex,epstopdf,url,algorithm,algorithmic,caption}
\usepackage{endnotes,authblk}
\usepackage{epsfig,subfigure,epstopdf}
\usepackage[colorlinks,linktocpage,linkcolor=blue]{hyperref}

\usepackage{etoolbox} 
\makeatletter
\patchcmd{\@makechapterhead}{\large}{\normalsize}{}{}
\patchcmd{\@makechapterhead}{\large}{\normalsize}{}{}
\patchcmd{\@makeschapterhead}{\normalsize}{\normalsize}{}{}
\makeatother
\usepackage[textwidth=6.0in,textheight=9in]{geometry}
\makeatletter
\g@addto@macro\normalsize{\setlength\abovedisplayskip{4pt}}
\g@addto@macro\normalsize{\setlength\belowdisplayskip{4pt}}
\makeatother

\newtheoremstyle{exampstyle}
  {\topsep} 
  {\topsep} 
  {\itshape} 
  {} 
  {\bfseries} 
  {.} 
  {.5em} 
  {} 

\theoremstyle{exampstyle}

\newtheorem{theorem}{Theorem}
\newtheorem{lemma}{Lemma}[section]
\newtheorem{corollary}{Corollary}[section]

\newtheorem{assumption}{Assumption}[section]
\newtheorem{remark}{Remark}[section]
\usepackage{parskip}
\let\oldref\ref
\renewcommand{\ref}[1]{(\oldref{#1})}  
\renewcommand{\eqref}[1]{(\oldref{#1})} 
\DeclareMathOperator*{\argmin}{arg\,min}

\newcommand{\blue}[1]{\textcolor{blue}{#1}}

\newcommand{\red}[1]{\textcolor{red}{#1}}

\newbox\boxaddrone \newbox\boxaddrtwo

\def\N+{n\in\mathbb{N}^{+}}

\def\n{\partial{\overrightarrow{\bf n}}}
\def\d{\mathcal{D}((-\Delta)^\gamma)}

\def\L{\mathcal{L}}

\def\l{\langle}
\def\ro{\rangle_{L^2(\Omega)}}

\def\S{\mathcal{S}}

\def\re{\operatorname{Re}}
\def\s{\text{Span}}
\begin{document}

\title{
On the identification of source term
in the heat equation from sparse data}
%
\author[1]{William Rundell\thanks{rundell@math.tamu.edu}}
\author[2]{Zhidong Zhang\thanks{zhidong.zhang@helsinki.fi}}
\affil[1]{Department of Mathematics, Texas A\&M University, USA}
\affil[2]{Department of Mathematics and Statistics, University of Helsinki, Finland}

\maketitle

\begin{abstract}
We consider the recovery of a source term $f(x,t)=p(x)q(t)$ for the
nonhomogeneous heat equation in $\Omega\times (0,\infty)$
where $\Omega$ is a bounded domain in $\mathbb{R}^2$ with smooth boundary
$\partial\Omega$ from overposed lateral data on a sparse subset of
$\partial\Omega\times(0,\infty)$.
Specifically, we shall require a small finite number $N$ of measurement
points on $\partial\Omega$ and prove a uniqueness result; namely
the recovery of the pair $(p,q)$ within a given class, by a judicious choice
of $N=2$ points.
Naturally, with this paucity of overposed data, the problem is severely
ill-posed. Nevertheless we shall show that provided the data noise level
is low, effective numerical reconstructions may be obtained.\\



\end{abstract}
\section{Introduction}

The inverse problem of recovering an unknown source term $f$ in the parabolic
equation $u_t - \triangle u = f$ from overspecified data on the solution $u$
has a long history, see for example,
\cite{Cannon:1968,Rundell:1983,Isakov:1990}.
A brief summary of this can be encapsulated by the observation that to obtain
a term $f=f(x,t)$ will require either knowledge of $u$ over
$\mathbb{R}^{n}\times\mathbb{R}$, which is impractical in almost every physical
situation, or over a sufficiently dense subset whereby an approximation
could be determined.
Thus most work has concentrated on one of the special cases
$f = q(t)$ or $f=p(x)$ or as a product $f = p(x)q(t)$ where either
$q$ or $p$ is known.
An exception here is \cite{CannonEsteva:1986} where the problem
was considered in $\mathbb{R}\times(0,\infty)$ and $p$ was of compact support.

It has also been observed that the recovery of a spatially unknown $f$
from spatial measurements of $u$ is usually only mildly-ill-conditioned
but the recovery from temporal measurements of $u$ is severely ill-posed.
The situation for $f=q(t)$ is reversed.
In fact, this problem spawned the now well-known notion that
to recover an unknown term or coefficient  in a partial differential equation
one should ideally prescribe data in a ``parallel'' direction to that of
the unknown; giving overposed data in the ``orthogonal'' direction
is likely to be severely ill-posed.

In this paper we shall assume the form $f = p(x)q(t)$ where both
$p$ and $q$ are unknown.  We shall prescribe extremely sparse 
time-trace data and show unique recovery within the specified spaces
in which $p$ and $q$ are defined although this paucity of data will require
quite severe restrictions on the allowable class for the  unknown term $q(t)$.
The exposition will be much simpler if we take the spatial domain $\Omega$
to be the unit disc in $\mathbb{R}^2$.
This is not an essential requirement and we could take
$\Omega\subset\mathbb{R}^2$ to have a smooth $C^2$ boundary.
We will comment on this fact later.

Let $u(x,t) \in \Omega\times[0,\infty)$ solve
\begin{equation}\label{eqn:direct_pde}
 \begin{cases}
  \begin{aligned}
   \frac{\partial u}{\partial t}(x,t)-\triangle u(x,t)&=p(x)q(t), &&(x,t)\in \Omega\times(0,\infty),\\
   u(x,t)&=0, &&(x,t)\in \partial\Omega\times(0,\infty),\\
   u(x,0)&=0, &&x\in \Omega.
  \end{aligned}
 \end{cases}
\end{equation}
As noted, $\Omega$ is the unit disc in $\mathbb{R}^2$ and $p$, $q$
are the unknown source sub-functions.
Our additional data is of the form of flux measurements at a small
number $L$ of points situated on $\partial\Omega$
\begin{equation}\label{eqn:overposed_data}
 g_\ell(t):=\frac{\partial u}{\n}(z_\ell,t),\ t\in (0,\infty),\ z_\ell\in\partial\Omega,
 \ \ell=1,2,\ldots\,,L. 
\end{equation} 

A related problem was considered in 
\cite{HettlichRundell:2001} where it was assumed that $q=1$ and $p=\chi(D)$
for some star-like domain $D\subset\Omega$.
Uniqueness in the form
of local injectivity of the derivative of the map $D\to g_\ell$
was shown: that is recovery of
the shape and location of a source of known uniform strength.
In this case only two flux measurements were required, that is $L=2$.

Our goal in this paper is to generalize this result to include
a nontrivial time-dependent term $q(t)$.
Such a modification represents a more realistic physical situation whereby
the strength of the source  may change with time.
We will show an analogous result to that in \cite{HettlichRundell:2001}
again requiring only that $L=2$ but it should not be surprising
that full generality cannot hold for $q(t)$.
We will show that uniqueness holds if $q(t)$ is a sequence of
step functions, that is $q(t)=\sum_{k=1}^K q_kH(t-c_k)$ where 
$H(t)$ is the Heaviside function
and $\{q_k,c_k\}_1^K$ will be determined in addition to $p(x)$. 
Note that the case of $K=\infty$ is allowed.
However, the previous result in \cite{HettlichRundell:2001}
has very little lee-way for generalization
and indeed we have unable to allow simply $p(x)\in L^\infty(\Omega)$ which
would open the path to allow $p=\chi(D)$ as before.
Instead we will have to assume slightly more regularity in $p$ and this
will preclude allowing it to represent the characteristic function
of a subset $D\subset\Omega$.
On the other hand we will be able to approximate this to within
any desired accuracy so from a physical standpoint there is no essential loss
of generality.
We shall show this case with some numerical runs in the final section.

There are many physical applications of this work and we mention only
the following.
Suppose there is an extended  source whose spatial location is only known
approximately. This could be a source of pollutant for example.
It is also a likely possibility that the output from this source depends on
time but that over a small enough period can be considered to be approximately
constant.
Measurements  can only be made at a distance from the source and the number
of measurement points is very small, perhaps due to logistics, but also
due to a small number of detecting sensors.
One could consider this problem to be in all of $\mathbb{R}^2$ or 
assume that it is more localized with given boundary constraints.
The latter situation is more complex and is the one taken in this paper.

Thus the main result of this paper is as follows, we shall describe some of the
technical definitions involved in the next section.
\begin{theorem}\label{thm:uniqueness}
Set the boundary observation points $\{z_\ell\}$ as 
$z_\ell=(\cos{\theta_\ell},\sin{\theta_\ell})$. 
 Then under Assumption \ref{assumption}, two boundary flux observations 
 can uniquely determine $(p,q)$ up to multiplication, provided  
 \begin{equation}\label{condition}
  \theta_1-\theta_2\notin \pi\mathbb{Q},
 \end{equation}
 where $\mathbb{Q}$ is the set of rational numbers. 
 
 More precisely, let $(p(x),q(t)),\ (\tilde{p}(x), \tilde{q}(t))$ satisfy Assumption 
 \ref{assumption} and denote the corresponding solutions as 
 $u, \tilde{u}$, respectively. If  
 $$\frac{\partial u}{\n}(z_\ell,t)=\frac{\partial \tilde{u}}{\n}(z_\ell,t),
 \quad t\in(0,\infty),\ \ell=1,2,$$  
 and the condition \eqref{condition} is fulfilled,  
 then there exists a constant $C_0\ne0$ such that 
 $p=C_0\tilde{p}$ in $L^2(\Omega)$ and $q=C_0^{-1}\tilde{q}$ on 
 $[0,\infty)$. 
\end{theorem}

This article is outlined as follows. In section 2, we provide several 
preliminary results and prove some lemmas which play crucial 
roles in the proof of the main theorem.
In section~3, we show the well-definedness 
and the analytic continuation of the Laplace transform on the flux data, 
see Lemma \ref{analytic}, and three auxiliary lemmas for the uniqueness 
proof. These allow the  completion of the proof of Theorem \ref{thm:uniqueness}
in section~3.4.
Based on the theoretical uniqueness property, we 
construct an iterative scheme to reconstruct the unknowns $p$, $q$ 
and several numerical results are reproduced in section~4.

\section{Preliminary lemmas and background}

\subsection{Eigensystem $\{\lambda_n, \varphi_n:\N+\}$}
Let $\{\lambda_n, \varphi_n(x):\N+\}$ be eigenpairs of $-\triangle$ on $\Omega$ 
with Dirichlet boundary conditions.
The corresponding eigenfunctions $\{\varphi_n\}$ will be used in
polar coordinates
\begin{equation}\label{eigenfunction}
\varphi_n(r,\theta)=\omega_nJ_{m(n)}(\sqrt{\lambda_n}r)\cos{(m(n)\theta+\phi_n)}.
\end{equation} 

\begin{remark}
 In the representation of $\varphi_n$, $\omega_n$ is the normalized 
 coefficient to make sure $\|\varphi_n\|_{L^2(\Omega)}=1$, $J_m$ is 
 the $m$-th order Bessel function and the phase $\phi_n$ is $0$ or $-\pi/2$. 
 The eigenvalues $\{\lambda_n\}$ are set as the square of the zeros 
 of Bessel functions $\{J_m\}$ with integer $m$. By Bourget's hypothesis, 
 which was proven in \cite{Siegel:2014}, there exists no common 
 positive zeros between two Bessel functions with integer orders. 
 After indexing all 
 the eigenvalues by non-decreasing order, with a fixed $n$, we can get 
 the corresponding value of $m$ so that we view  $m$ as a function of $n$ 
 and this dependence is reflected in the notation $m(n)$. 
 Since there are two choices $0$ or $-\pi/2$ for $\phi_n$, for each eigenvalue 
 $\lambda_n$ with nonzero $m(n)$, the multiplicity is two.
Thus it has two corresponding eigenfunctions.
When $m(n)=0$, the multiplicity is only one since the angular part 
 $\cos{(m(n)\theta-\pi/2)}=\sin{0}$ vanishes on $[0,2\pi)$. 
This fact is also guaranteed in the more general case of a non-circular
domain $\Omega$ by the Krein-Rutman theorem.
 For more details on the structure of the eigenfunctions, see 
 \cite{GrebenkovNguyen:2013}.
 
 Since $-\Delta$ is self-adjoint and positive definite on $\Omega$  
 with homogeneous Dirichlet boundary condition, $\{\lambda_n\}$ will 
 be strictly positive and $\{\varphi_n\}$ constitutes an orthonormal 
 basis of $L^2(\Omega)$.
We index the eigenvalues with non-decreasing order, then we have 
 \begin{equation*}
  0<\lambda_1\le\cdots\le \lambda_n\le \cdots\ \text{(multiplicity counted)}, 
  \quad \lim_{n\to \infty} \lambda_n=\infty.
 \end{equation*}  
\end{remark}

\begin{remark}\label{remark:lambda&bessel}
 Here we list some properties of $\{\lambda_n\}$ and $J_m(x)$ which 
 will be used later. 
 \begin{itemize}
  \item By Weyl's law, $\Omega\subset\mathbb{R}^2$ implies 
  $\lambda_n=O(n)$.
  \item $[x^{m+1}J_{m+1}(x)]'=x^{m+1}J_m(x)$.
  \item $2mJ_m(x)/x=J_{m-1}(x)+J_{m+1}(x)$.
  \item $2[J_m(x)]'=J_{m-1}(x)-J_{m+1}(x)$. 
 \end{itemize}
\end{remark}

The following lemma concerns the estimate for the normalized coefficient 
$\omega_n$.
\begin{lemma}\label{omega_n}
For $\N+,$ $\{\omega_n\}$ are given by
 \begin{equation*}
 \omega_n = 
 \begin{cases}
  2^{1/2}\pi^{-1/2}[J_{m+1}(\lambda_n^{1/2})]^{-1}, & m(n)\ne 0,\\
  \pi^{-1/2}[J_1(\lambda_n^{1/2})]^{-1}, & m(n)=0. 
 \end{cases}
\end{equation*}
\end{lemma}
\begin{proof}
 Since $\|\varphi_n\|_{L^2(\Omega)}=1,$ then 
 \begin{equation*}
  \omega_n^2 \int_0^{2\pi} \cos^2{(m\theta+\phi_n)}\,d\theta 
  \ \int_0^1 J_m^2(\sqrt{\lambda_n}r)r\,dr=1.
 \end{equation*}
If $m(n)\ne 0$, we have 
$ \int_0^{2\pi} \cos^2{(m\theta+\phi_n)}\,d\theta
 =\frac{1}{2}\int_0^{2\pi} {1+\cos{2(m\theta+\phi_n)}}\,d\theta=\pi$
and 
\begin{equation*}
\begin{aligned}
 \int_0^1 J_m^2(\sqrt{\lambda_n}r)r\,dr&=\lambda_n^{-1}\int_0^{\lambda_n^{1/2}} 
 J_m^2(s)s \,ds \\
 &=\lambda_n^{-1}\Big[r^2J_{m+1}^2(r)/2+r^2J_m^2(r)/2
 -mrJ_m(r)J_{m+1}(r)\Big]\Big|_0^{\lambda_n^{1/2}}\\
 &=J_{m+1}^2(\lambda_n^{1/2})/2,
 \end{aligned}
\end{equation*}
where the second result comes from the fact that $\lambda_n^{1/2}$ is the 
zero of $J_m(r)$ and the recurrence relations in Remark 
\ref{remark:lambda&bessel}. 
Hence, we have 
\begin{equation*}
 \omega_n J_{m+1}(\lambda_n^{1/2})= 2^{1/2}\pi^{-1/2}. 
\end{equation*}
Analogously, for the case of $m(n)=0$, it holds that 
$\omega_n J_1(\lambda_n^{1/2})= \pi^{-1/2}$ and the proof is complete. 
\end{proof}

\subsection{Assumptions and solution regularities}

We give the definitions of the space $\d$ and the Heaviside 
function $H(t)$ that will be used throughout the paper.
For $\gamma>0$, define $\d\subset L^2(\Omega)$ as 
\begin{equation}\label{eqn:D_gamma}
 \d:=\{\psi\in L^2(\Omega): \sum_{n=1}^\infty \lambda_n^{2\gamma} 
 \ |\l\psi(\cdot),\varphi_n(\cdot)\ro|^2\}<\infty,
\end{equation}
here $\l\cdot,\cdot\ro$ means the inner product in $L^2(\Omega)$. 
Due to $\Omega$ is the unit disc, $\d\subset H^{2\gamma}(\Omega)$. 
Given $0<\gamma_1<\gamma_2$, since $0<\lambda_1\le \cdots\le \lambda_n\le \cdots,$  
it is not hard to show 
$\mathcal{D}((-\Delta)^{\gamma_2})\subset\mathcal{D}((-\Delta)^{\gamma_1}).$ 

Also, the Heaviside function $H(t)$ is defined in the usual way
\begin{equation*}
 H(t)=
 \begin{cases}
  0,&t<0,\\
  1,&t\ge0,
 \end{cases}
\end{equation*}
and it is clear that $\chi_{{}_{[a,b)}}=H(t-a)-H(t-b),\ a<b.$

With these definitions, we require
the following assumptions to be valid throughout the paper.  
\begin{assumption}\label{assumption}
$p(x)$ and $q(t)$ satisfy the following conditions: 
\begin{itemize}
\item There exists $\gamma>0$ such that $p(x)\in \d$ and
$\|p\|_{L^2(\Omega)}\ne 0$.
  \item $q(t)\in L^1(0,\infty)$ is a piecewise constant function, 
  or written as a linear combination of Heaviside functions, 
  \begin{equation*}
  q(t)=\sum_{k=1}^{K-1}\beta_k \chi_{{}_{[c_k,c_{k+1})}}=\sum_{k=1}^K q_kH(t-c_k),
  \end{equation*}
  where 
  \begin{equation*}
   K\in\mathbb{N}^+\cup\{\infty\},\quad 0\le c_1<c_2<\cdots, \quad q_k\ne 0.
  \end{equation*}
  Moreover, there exists $\eta>0$ such that 
  \begin{equation}\label{eta}
   \inf\{|c_k-c_{k+1}|:k=1,\cdots,K-1\}\ge \eta.
  \end{equation}  
\end{itemize}
\end{assumption}

\begin{remark}\label{condition_q}
 The inclusion $q(t)\in L^1(0,\infty)$ and the infimum $\eta$ give that 
 \begin{equation*}
   \eta\sum_{k=1}^{K-1} |\beta_k|\le \sum_{k=1}^{K-1} |\beta_k|\ |c_k-c_{k+1}| 
   = \|q\|_{L^1(0,\infty)}< \infty,
 \end{equation*}
 which leads to $\sum_{k=1}^{K-1} |\beta_k|<\infty$. 
 Also, we have 
 \begin{equation*}
 \begin{aligned}
  \|q\|^2_{L^2(0,\infty)}=&\sum_{k=1}^{K-1} |\beta_k|^2\ |c_k-c_{k+1}|
  =\sum_{k=1}^{K-1} (|\beta_k|\ |c_k-c_{k+1}|)\ (|\beta_k|)\\
  \le &\Big[\sum_{k=1}^{K-1} |\beta_k|\ |c_k-c_{k+1}|\Big]
  \ \ \Big[\sum_{k=1}^{K-1} |\beta_k|\Big]<\infty,
  \end{aligned}
 \end{equation*}
which gives $q\in L^2(0,\infty)$.
 
 From the equality $\sum_{k=1}^{K-1}\beta_k \chi_{{}_{[c_k,c_{k+1})}}
 =\sum_{k=1}^K q_kH(t-c_k)$, we derive that for $K<\infty$, 
 \begin{equation*}
  q_1=\beta_1,\ q_K=-\beta_{K-1}, 
  \ q_k=\beta_k-\beta_{k-1}\ \ \text{for}\ \  2\le k\le K-1, 
 \end{equation*}
 and for $K=\infty$, 
  \begin{equation*}
  q_1=\beta_1, \ q_k=\beta_k-\beta_{k-1}\ \ \text{for}\ \  k\ge 2. 
 \end{equation*}
Hence, 
 \begin{equation*}
  |q_k|\le \sum_{k=1}^{K-1} |\beta_k|<\infty,
  \quad \sum_{k=1}^K|q_k|\le 2 \sum_{k=1}^{K-1} |\beta_k|<\infty.
 \end{equation*}

In addition, \eqref{eta} yields that 
\begin{equation*}
 |c_{k_1}-c_{k_2}|\ge |k_1-k_2|\eta, \quad k_1,k_2=1,\cdots,K. 
\end{equation*}
\end{remark}

In this subsection we also give a regularity result for 
$\frac{\partial u}{\n}(\cdot,t)$.

\begin{lemma}\label{lemma:flux regularity}
For a.e. $t\in[0,\infty)$, 
	$\frac{\partial u}{\n}(\cdot,t)\in C^{0,2\gamma}(\partial\Omega).$
\end{lemma}
\begin{proof}
 From Assumption \ref{assumption}, we have $p(x)\in \d$ and $\gamma>0$. 
 Since $\mathcal{D}((-\Delta)^{\gamma_2})\subset\mathcal{D}((-\Delta)^{\gamma_1})$ 
 if $0<\gamma_1<\gamma_2$, we can set $\gamma\in (0,1/4)$. In addition, 
 recalling $q(t)\in L^2(0,\infty)$, then we have 
 $p(x)q(t)\in L^2(0,\infty;\d)$. From the spectral representation of $u$, 
 the following regularity holds
 $$u\in L^2(0,\infty;\mathcal{D}((-\Delta)^{\gamma+1}))
 \subset L^2(0,\infty;H^{2\gamma+2}(\Omega)).$$
 Then using the continuity of the trace map 
 $\psi\in H^{2\gamma+2}(\Omega)\mapsto \frac{\partial \psi}{\n}\in 
 H^{2\gamma+1/2}(\partial\Omega),$ which is \cite[Theorem 9.4]{LionsMagenes:1972V1},  
 gives that for a.e. $t\in [0,\infty)$, 
 $\frac{\partial u}{\n}(\cdot,t)\in H^{2\gamma+1/2}(\partial\Omega)$. 
 Note that $\partial\Omega$ is one-dimensional and 
 $2\gamma+1/2\in(1/2,1)$, which mean the conditions of  
 \cite[Theorem 8.2]{DiNezzaPalatucciValdinoci:2012} are satisfied. 
 Then we have $\frac{\partial u}{\n}(\cdot,t)\in 
 C^{0,2\gamma}(\partial\Omega)$ and this completes the proof. 
\end{proof}


\section{Uniqueness}
This section is devoted to the proof of the main theoretical result,
Theorem \ref{thm:uniqueness}.

\subsection{Harmonic functions and measurements representations}
First, we need to show how to connect the boundary flux 
measurements $\frac{\partial u}{\n}|_{\partial\Omega}$ and the unknowns $p(x),\ q(t)$.  
Here we introduce the harmonic functions $\{\xi_j:j\in\mathbb{N}^+\}$ 
which will be used to represent measurements. The set of harmonic functions 
with domain $(r,\theta)\in[0,1]\times[0,2\pi)$ is defined as  
\begin{equation*}
\begin{aligned}
 \xi_j(r,\theta)&=
 \begin{cases}
 \pi^{-1/2}r^l\cos{l\theta}, & j=2l+1,\ l>0,\\ 
 2^{-1/2}\pi^{-1/2}, & j=1,\\
  \pi^{-1/2}r^l\sin{l\theta},& j=2l,\ l>0,  
 \end{cases}\\
 &=
 \begin{cases}
  \pi^{-1/2}r^{\lfloor j/2\rfloor}\cos{\big(\lfloor j/2\rfloor\theta
 +\sigma_j\big)}, &j>1,\\
 2^{-1/2}\pi^{-1/2}, & j=1.
 \end{cases}
 \end{aligned}
\end{equation*}
Here $\lfloor j/2 \rfloor$
means the largest integer which is not larger than $j/2$ and
\begin{equation*}
 \sigma_j= \begin{cases}
 0, &j\ \text{is odd},\\ 
 -\pi/2, & j\ \text{is even}.
 \end{cases}
\end{equation*}
Fixing $r=1,$ the set 
$$\{\xi_j(1,\theta): j\in\mathbb{N}^+,\ \theta\in[0,2\pi)\}
=\{2^{-1/2}\pi^{-1/2},\ \pi^{-1/2}\sin{(l\theta)},
\ \pi^{-1/2}\cos{(l\theta)}:l\in\mathbb{N}^+\}$$ 
will form an orthonormal basis in $L^2(\partial \Omega)$.  

Fix $z=(\cos{\theta_z},\sin{\theta_z})\in \partial\Omega$, let's calculate 
$\l \xi_j(z)\xi_j(\cdot), \varphi_n(\cdot)\ro$. For $\N+$ define 
 \begin{equation}\label{a_np_n} 
 \begin{aligned} 
 a_n(z)&=
 \begin{cases}
  2^{1/2}\pi^{-1/2}\lambda_n^{-1/2} \cos{(m(n)\theta_z+\phi_n)},&m(n)\ne0,\\
  \pi^{-1/2}\lambda_n^{-1/2},  &m(n)= 0,
 \end{cases}\\
  p_n&=\l p(\cdot),\varphi_n(\cdot)\ro.
 \end{aligned}
 \end{equation}
Reference~\cite{HettlichRundell:2001} shows that 
$\l \xi_j(z)\xi_j(\cdot), \varphi_n(\cdot)\ro $ can be written as a product of
three terms:
$$
\cos{(\lfloor j/2\rfloor\theta_z+\sigma_j)},\quad
\int_0^{2\pi}\!\cos(m\theta + \phi_n)\cos{(\lfloor j/2\rfloor\theta+\sigma_j)}\,d\theta,
\quad 
\omega_n\int_0^1 r^{\lfloor j/2\rfloor+1} J_m(\lambda_n^{1/2}r)\,dr
$$
and a factor of $\pi^{-1}$ or $2^{-1}\pi^{-1}$. 
The integral in the angular variable $\theta$, and hence the inner product,
is zero except when $m=\lfloor j/2\rfloor$
and $\phi_n=\sigma_j$ in which case it has value $\pi$ or $2\pi$. 
The integral in the radial variable $r$ can then be written as 
$\int_0^1 r^{m+1} J_m(\lambda_n^{1/2}r)\,dr$
and after a change of variable $s=\lambda_n^{1/2}r$ and use of a Bessel
function recursion formula, becomes
$$\lambda_n^{-1-m/2}\int_0^{\lambda_n^{1/2}} s^{m+1}J_m(s)\,ds =
\lambda_n^{-1/2}J_{m+1}(\lambda_n^{1/2}).$$

Thus combining all the terms
and using Lemma \ref{omega_n} show that 
 \begin{equation}\label{a}
 \begin{aligned}
 \l \xi_j(z)\xi_j(\cdot), \varphi_n(\cdot)\ro
 =&\begin{cases}
   (\frac{2}{\pi\lambda_n})^{1/2} \cos{(m\theta_z+\phi_n)},  
   &\text{if}\ m(n)= \lfloor j/2\rfloor\ne0\ \text{and}\ \phi_n=\sigma_j,\\
 (\pi\lambda_n)^{-1/2},  &\text{if}\ m(n)= \lfloor j/2\rfloor=0,\\
 0, &\text{otherwise},
   \end{cases}\\
 =&\begin{cases}
   a_n(z),  &\text{if}\ m(n)= \lfloor j/2\rfloor\ \text{and}\ \phi_n=\sigma_j,\\
 0, &\text{otherwise}.
   \end{cases}
 \end{aligned}
 \end{equation}

Now we use the Harmonic basis $\{\xi_j:j\in\mathbb{N}^+\}$ 
to build a connection between the boundary flux $\partial u(z,t)/\n$ 
and source terms $p(x),q(t)$. 

Fix a point $z\in \partial\Omega$  and define 
$\psi_z^M\in C^\infty(\overline{\Omega})$ as  
\begin{equation}\label{eqn:psi_definition}
 \psi_z^M(x)=\sum_{j=1}^M \xi_j(z)\xi_j(r,\theta),
 \quad (r,\theta)\in [0,1]\times[0,2\pi).
\end{equation}
Then we denote the solutions of the following systems by $\{u_z^M\}$, 
\begin{equation}\label{eqn:u_singular}
\begin{cases}
 \begin{aligned}  
  \frac{\partial u_z^M}{\partial t}(x,t)-\triangle u_z^M(x,t)&=0,
  &&(x,t)\in \Omega\times(0,\infty),\\
 u_z^M(x,t)&=0, &&(x,t)\in \partial\Omega\times(0,\infty),\\
 u_z^M (x,0)&=-\psi_z^M, &&x\in \Omega,
 \end{aligned}
\end{cases}
\end{equation}
and require the lemma below. 

\begin{lemma}\label{data_1}
 Let $w_z^M=u_z^M+\psi_z^M,$ then we have 
 \begin{equation*}
 -\int_0^t \frac{\partial u}{\n}(z, \tau)\ d\tau
 =\int_0^t q(\tau)\ \Big[\lim_{M\to \infty}\big\l p(\cdot), 
 w_z^M(\cdot,t-\tau)\big\rangle_{L^2(\Omega)}\Big]\ d\tau.
\end{equation*}
\end{lemma}
\begin{proof}
 Since $\psi_z^M$ is the linear combination of harmonic functions on $\Omega,$ then 
 \begin{equation*}
  \frac{\partial\psi_z^M}{\partial t} -\triangle \psi_z^M=0.  
 \end{equation*}
This result and \eqref{eqn:u_singular} show that $w_z^M$ satisfies the equation 
\begin{equation*}
 \frac{\partial w_z^M}{\partial t} -\triangle w_z^M=0, 
 \ \  (x,t)\in \Omega\times(0,\infty)
\end{equation*}
with zero initial condition and the boundary condition 
$w_z^M|_{\partial\Omega}=\psi_z^M|_{\partial\Omega}$. Thus Green's identities show that 
for each $v\in H_0^1(\Omega)$, 
\begin{equation*}
 \int_\Omega \frac{\partial w_z^M}{\partial t}(x,t)\ v(x)
 +\triangledown w_z^M(x,t) \cdot \triangledown v(x)\ dx=0,\ \ t\in (0,\infty).
\end{equation*}

By a direct calculation we obtain 
\begin{equation*}
 \begin{aligned}
  \int_0^t \int_\Omega p(x) q(\tau) w_z^M(x,t-\tau)\ dx\ d\tau
  =\int_0^t\int_\Omega \left[\frac{\partial u}{\partial t}(x,\tau)
  -\triangle u(x,\tau)\right]\ w_z^M(x,t-\tau)\ dx\ d\tau. 
 \end{aligned}
\end{equation*}
Green's identities and the vanishing initial conditions of $u$ and 
$w_z^M$ give that 
\begin{equation*}
 \begin{aligned}
 \int_0^t\int_\Omega \frac{\partial u}{\partial t}(x,\tau)\ w_z^M(x,t-\tau)\ dx\ d\tau
 =& \int_0^t\int_\Omega \frac{\partial w_z^M}{\partial t} (x,t-\tau)u(x,\tau)\ dx\ d\tau,\\
 \int_0^t\int_\Omega -\triangle u(x,\tau)\ w_z^M(x,t-\tau)\ dx\ d\tau
 =& \int_0^t\int_\Omega \triangledown u(x,\tau)\cdot\triangledown 
 w_z^M(x,t-\tau)\ dx\ d\tau \\
 &-\int_0^t \int_{\partial\Omega} \frac{\partial u}{\n}(x,\tau)
 \psi_z^M(x)\ dx\ d\tau.
\end{aligned}
\end{equation*}
Hence, 
\begin{equation*}
 \begin{aligned}
  \int_0^t \int_\Omega p(x) q(\tau) w_z^M(x,t-\tau)\ dx\ d\tau
  =&\int_0^t\int_\Omega \Big[\frac{\partial w_z^M}{\partial t} (x,t-\tau)\ u(x,\tau)\\
  &\qquad\quad+\triangledown w_z^M(x,t-\tau)\cdot\triangledown u(x,\tau)\Big]
  \ dx\ d\tau\\
 & -\int_0^t \int_{\partial\Omega} \frac{\partial u}{\n}(x,\tau)
 \psi_z^M(x)\ dx\ d\tau\\
 =&-\int_0^t \sum_{j=1}^M \xi_j(z)\ \big\l\frac{\partial u}{\n}(\cdot,\tau),\xi_j(\cdot) 
 \big\rangle_{L^2(\partial\Omega)} \ d\tau.
 \end{aligned}
\end{equation*}
The smoothness property  $\frac{\partial u}{\n}(\cdot,t)\in 
C^{0,2\gamma}(\partial\Omega)$ in Lemma \ref{lemma:flux regularity} 
ensures that the Fourier series of $\frac{\partial u}{\n}(\cdot,t)$ 
converges pointwisely on $\partial\Omega$, namely 
\begin{equation*}
 \lim_{M\to \infty}\sum_{j=1}^M \xi_j(z)\ \big\l\frac{\partial u}{\n}(\cdot,\tau),\xi_j(\cdot) 
 \big\rangle_{L^2(\partial\Omega)} =\frac{\partial u}{\n}(z,\tau),\quad a.e.\ \tau\in (0,t).
\end{equation*}
Since the `almost everywhere' does not effect the result of integral, 
we have  
\begin{equation*}
\begin{aligned}
 -\int_0^t \frac{\partial u}{\n}(z, \tau)\ d\tau
 &=-\int_0^t \Big[\lim_{M\to \infty}\sum_{j=1}^M \xi_j(z)\ \big\l\frac{\partial u}{\n}(\cdot,\tau),\xi_j(\cdot) 
 \big\rangle_{L^2(\partial\Omega)}\Big]\ d\tau\\
 &=\int_0^t q(\tau)\ \Big[\lim_{M\to \infty}\big\l p(\cdot), 
 w_z^M(\cdot,t-\tau)\big \rangle_{L^2(\Omega)}\Big]\ d\tau.
 \end{aligned}
\end{equation*}
\end{proof}

With the above lemma, the next corollary follows.

\begin{corollary}\label{data_formula}
Fix $z=(\cos{\theta_z},\sin{\theta_z})\in \partial\Omega$, then 
\begin{equation*}
 -\int_0^t \frac{\partial u}{\n}(z, \tau)\ d\tau=\int_0^t q(\tau)
 \ \Big[\sum_{n=1}^{\infty} a_n(z)p_n [1-e^{-\lambda_n (t-\tau)}]\Big]\ d\tau.
\end{equation*}
\end{corollary}

\begin{proof}
 For each $M\in \mathbb{N}^+$, from \eqref{eqn:psi_definition} we have 
 $\psi_z^M\in L^2(\Omega)$. Then the Fourier expansion of 
 $\psi_z^M$ can be given as $\psi_z^M=\sum_{n=1}^\infty a_n^M(z) \varphi_n$, and from 
 $w_z^M=u_z^M+\psi_z^M$ we obtain 
 \begin{equation*}
  w_z^M(x,t)=\sum_{n=1}^\infty a_n^M(z) (1-e^{-\lambda_n t})\varphi_n(x).
 \end{equation*}
 The above representation and the regularity $\psi_z^M\in L^2(\Omega)$ 
 give that $w_z^M(\cdot,t)\in L^2(\Omega)$ for $t\in [0,\infty)$. 
 Since $p, w_z^M(x,t)$ both belong to $L^2(\Omega)$, we have 
 \begin{equation*}
  \big\l p(\cdot), w_z^M(\cdot,t)\big \rangle_{L^2(\Omega)}
  =\sum_{n=1}^\infty a_n^M(z)p_n (1-e^{-\lambda_n t}).
 \end{equation*}
 \eqref{a} shows that 
 \begin{equation*}
  a_n^M(z)=
  \begin{cases}
   a_n(z), &m(n)< M/2,\\
   0,&m(n)>M/2, 
  \end{cases}
 \end{equation*}
 and for the case of $m(n)=M/2$, 
  \begin{equation*}
  a_n^M(z)=
  \begin{cases}
   a_n(z), &\sigma_n=-\pi/2,\\
   0,&\sigma_n=0. 
  \end{cases}
 \end{equation*}
 These results mean that $a_n^M(z)=a_n(z)$ if $M$ is large, and 
 $|a_n^M(z)|\le |a_n(z)|$ for each $n,M$. 
 
 Given $\epsilon>0$, Lemma \ref{convergence}, which will be proved in the 
 next subsection, yields that there 
 exists large $l>0$ such that $\sum_{n=l}^\infty |a_n(z)p_n|<\epsilon$. 
 From the above results for $a_n^M(z)$, we can find a $M_0$ such that 
 if $M\ge M_0$, $a_n^M(z)=a_n(z)$ for $n=1,\cdots,l-1$. So for $M\ge M_0$, 
 \begin{equation*}
 \begin{aligned}
  \Big|\sum_{n=1}^\infty a_n^M(z)p_n (1-e^{-\lambda_n t})
  -\sum_{n=1}^\infty a_n(z)p_n (1-e^{-\lambda_n t})\Big|
  =&\Big|\sum_{n=l}^\infty [a_n^M(z)-a_n(z)]p_n (1-e^{-\lambda_n t})\Big|\\
  \le& 2\sum_{n=l}^\infty [|a_n^M(z)|+|a_n(z)|]\ |p_n| \\
  \le& 4\sum_{n=l}^\infty |a_n(z)p_n|<4\epsilon.
  \end{aligned}
 \end{equation*}
This in turn leads to 
\begin{equation*}
  \lim_{M\to \infty}\big\l p(\cdot), w_z^M(\cdot,t)\big \rangle_{L^2(\Omega)}
  =\lim_{M\to \infty}\sum_{n=1}^\infty a_n^M(z)p_n (1-e^{-\lambda_n t})
  =\sum_{n=1}^\infty a_n(z)p_n (1-e^{-\lambda_n t}),
\end{equation*}
which together with Lemma \ref{data_1} completes the proof. 
\end{proof}

\subsection{A Laplace transform analysis}
The uniqueness proof relies on the Laplace transform on the result 
in Corollary \ref{data_formula}. Before to analyze the Laplace transform, 
we need the following absolute convergence result. 
\begin{lemma}\label{convergence}
 $\sum_{n=1}^\infty a_n(z)p_n$ is absolute convergent for each $z\in\partial\Omega$.
\end{lemma}

\begin{proof}
By the Cauchy-Schwartz inequality, 
\begin{equation*}
\begin{aligned}
 \sum_{n=1}^\infty |a_n(z)p_n| 
 &=\sum_{n=1}^\infty |a_n(z)\lambda_n^{-\gamma}|
 \ |\lambda_n^{\gamma}p_n| \\
 &\le \Big[\sum_{n=1}^\infty a^2_n(z)\lambda_n^{-2\gamma}\Big]^{1/2}
 \ \Big[\sum_{n=1}^\infty \lambda_n^{2\gamma}p_n^2\Big]^{1/2}.
 \end{aligned}
\end{equation*}

$p(x)\in \d$ means that $\sum_{n=1}^\infty 
\lambda_n^{2\gamma}p_n^2<\infty.$ Also, from $\lambda_n=O(n)$, we have   
 \begin{equation*}
  a^2_n(z)\lambda_n^{-2\gamma} \le C|\lambda_n^{-1/2} 
  \lambda_n^{-\gamma}|^2
  \le C n^{-1-2\gamma},
 \end{equation*}
which gives $\sum_{n=1}^\infty a^2_n(z)\lambda_n^{-2\gamma}<\infty$. 
Hence we have $\sum_{n=1}^\infty |a_n(z)p_n|<\infty$ and completes the proof.
\end{proof}

Then from Corollary \ref{data_formula}, taking Laplace transform on 
$-\int_0^t \frac{\partial u}{\n}(z, \tau)\,d\tau$ w.r.t $t$ gives 
that  
\begin{equation*}
 \L \big(-\int_0^t \frac{\partial u}{\n}(z, \tau)\,d\tau\big)(s)
 =\L\big(q(t)\big)(s) \ \  \L\big(\sum_{n=1}^{\infty} a_n(z)p_n 
 [1-e^{-\lambda_n t}]\big)(s).
\end{equation*}
Since $q(t)\in L^1(0,\infty)$ is a piecewise constant function, 
it is bounded and $\L\big(q(t)\big)(s)$ 
is convergent and well-defined for $\re s>0$.
Again, it follows directly that
\begin{equation*}
 s\L\big(q(t)\big)(s)=\sum_{k=1}^K q_k e^{-c_k s},\ \re s>0.
\end{equation*}
From Lemma \ref{convergence} and $|1-e^{-\lambda_nt}|\le 1$, the series 
$\sum_{n=1}^{\infty} a_n(z)p_n [1-e^{-\lambda_n t}]$ is also uniformly bounded 
on $(0,\infty)$. This means its Laplace transform is well-defined 
for $\re s>0$ and the dominated convergence theorem can be applied to calculate 
the transform as  
\begin{equation}\label{laplace_analysis}
\begin{aligned}
 \L\big(\sum_{n=1}^{\infty} a_n(z)p_n[1-e^{-\lambda_n t}]\big)(s)
 &=\int_0^\infty \sum_{n=1}^{\infty} a_n(z)p_n[e^{-st}-e^{-(s+\lambda_n) t}]\,dt\\
 &=\sum_{n=1}^{\infty} a_n(z)p_n\int_0^\infty e^{-st}-e^{-(s+\lambda_n) t}\,dt\\
 &=\sum_{n=1}^{\infty} a_n(z)p_n\lambda_ns^{-1}(s+\lambda_n)^{-1},\quad \re s>0.
 \end{aligned}
\end{equation}
Now we have 
\begin{equation}\label{laplace}
 s^2\ \L \big(-\int_0^t \frac{\partial u}{\n}(z, \tau)\,d\tau\big)(s)
 =\Big[\sum_{k=1}^K q_k e^{-c_k s}\Big]\ 
 \Big[\sum_{n=1}^{\infty} a_n(z)p_n\lambda_n(s+\lambda_n)^{-1}\Big],\quad \re s>0.
\end{equation}
We will show the well-definedness and the analyticity for the above
complex-valued functions. 

\begin{lemma}\label{analytic}
 Under Assumption \ref{assumption}, the following properties hold. 
 \begin{itemize}   
 \item [(a)]  For $R\in\mathbb{R}$, define 
 $\mathbb{C}_R:=\{s\in \mathbb{C}:\re s>R\}$. Then 
 $\sum_{n=1}^\infty a_n(z)p_n\lambda_n(s+\lambda_n)^{-1}$ is uniformly convergent 
 for $s\in \mathbb{C}_R\setminus\{-\lambda_n:\N+\}$. 
 \item [(b)] $\sum_{n=1}^\infty a_n(z)p_n\lambda_n(s+\lambda_n)^{-1}$ is analytic 
 on $\mathbb{C}\setminus\{-\lambda_n:\N+\}$. 
 \item [(c)] $\sum_{k=1}^K q_k e^{-c_k s}$ is analytic on 
 $\mathbb{C}^+:=\{ s\in \mathbb{C}:\re s\ge0\}$.
 \end{itemize}
\end{lemma}
\begin{proof}
For $(a)$, since $0<\lambda_1\le\cdots\le\lambda_n\le\cdots\to\infty$, 
there exists a large $N_1$ such that $\lambda_n>2|R|$ for 
$n\ge N_1$. Then for $s\in \mathbb{C}_R\setminus\{-\lambda_n:\N+\}$ 
and $n\ge N_1$, 
\begin{equation*}
 |s+\lambda_n|\ge |\re s+\lambda_n|=\lambda_n+\re s\ge \lambda_n-|R|>0,
\end{equation*}
which gives 
\begin{equation*}
 |\lambda_n(s+\lambda_n)^{-1}|=\lambda_n|s+\lambda_n|^{-1}
 \le \lambda_n(\lambda_n-|R|)^{-1}< 2.
\end{equation*}
Given $\epsilon>0,$ Lemma \ref{convergence} yields that there exists $N_2$ such that 
for $l\ge N_2$, $$\sum_{n=l}^\infty|a_n(z)p_n|<\epsilon.$$
So, for $l\ge \max\{N_1,N_2\}$ and $s\in \mathbb{C}_R\setminus\{-\lambda_n:\N+\}$, 
\begin{equation*}
 \big|\sum_{n=l}^\infty a_n(z)p_n\lambda_n(s+\lambda_n)^{-1}\big| 
 \le \sum_{n=l}^\infty |a_n(z)p_n|\ |\lambda_n(s+\lambda_n)^{-1}|
 \le 2\sum_{n=l}^\infty |a_n(z)p_n|<2\epsilon,
\end{equation*}
which implies the uniform convergence. \\

For $(b)$, it is clear that $a_n(z)p_n\lambda_n(s+\lambda_n)^{-1}$ 
is holomorphic on $\mathbb{C}_R$. Then the uniform convergence gives that 
$\sum_{n=1}^\infty a_n(z)p_n\lambda_n(s+\lambda_n)^{-1}$ is holomorphic, 
i.e. analytic on  $\mathbb{C}_R$ for each $R\in\mathbb{R}$. Given 
$s\in\mathbb{C}\setminus\{-\lambda_n:\N+\}$, we can find a $R$ such that 
$s\in \mathbb{C}_R$, which means 
$\sum_{n=1}^\infty a_n(z)p_n\lambda_n(s+\lambda_n)^{-1}$ is analytic on 
$\mathbb{C}\setminus\{-\lambda_n:\N+\}$.\\

For $(c)$, it is obviously valid if $K<\infty$. This is because 
 $q_k e^{-c_k s}$ is analytic on $\mathbb{C}^+$ and the sum is finite.
For the case of $K=\infty$, following the proofs for $(a)$ and $(b)$, 
we have 
\begin{equation*}
 \big|\sum_{k=l}^\infty q_k e^{-c_k s}\big|\le \sum_{k=l}^\infty |q_k|\ |e^{-c_k s}| 
 \le \sum_{k=l}^\infty |q_k|,\quad s\in \mathbb{C}^+,\ l\in \mathbb{N}^+.
\end{equation*}
This result together with the absolute convergence of $\sum_{k=1}^\infty q_k$, 
stated by Remark \ref{condition_q}, yields the uniform convergence 
of $\sum_{k=1}^\infty q_k e^{-c_k s}$ on $\mathbb{C}^+$. Then with 
the analyticity of each component function $q_k e^{-c_k s}$, we can 
deduce that $\sum_{k=1}^\infty q_k e^{-c_k s}$ is analytic on 
$\mathbb{C}^+$ and complete the proof.
\end{proof}

\subsection{Auxiliary lemmas}
In order to prove Theorem \ref{thm:uniqueness}, 
some auxiliary lemmas are needed and stated below. 
\begin{lemma}\label{lemma_uniqueness_1}
 Write $z_\ell$ as $z_\ell=(\cos{\theta_\ell},\sin{\theta_\ell}),
 \ \ell=1,2$ and denote the set of distinct eigenvalues with 
 increasing order by $\{\lambda_j:j\in\mathbb{N}^+\}$.
 Provided the condition
 $\theta_1-\theta_2\notin \pi\mathbb{Q},$ where $\mathbb{Q}$ is 
 the set of rational numbers, then 
 $$\sum_{\lambda_n=\lambda_j} a_n(z_\ell)p_n=0,
 \quad j\in\mathbb{N}^+,\ \ell=1,2$$ implies that $p_n=0$ for $\N+$.
\end{lemma}
\begin{proof}
Fix $j\in \mathbb{N}^+$, if $m(n(j))\ne0$, then
 \begin{equation*}
 \sum_{\lambda_n=\lambda_j} a_n(z_\ell)p_n
 = 2^{1/2}\pi^{-1/2}\lambda_j^{-1/2}\bigl(\cos(m\theta_\ell)p_{n(j)}
 +\sin(m\theta_\ell)p_{n(j)+1}\bigr) =0,\quad \ell=1,2.
\end{equation*}
This means 
\begin{equation*}
 \begin{bmatrix}
  \cos(m\theta_1)&\sin(m\theta_1)\\
  \cos(m\theta_2)&\sin(m\theta_2)
 \end{bmatrix}
\begin{bmatrix}
 p_{n(j)}\\
 p_{n(j)+1}
\end{bmatrix}
=\begin{bmatrix}
  0\\0
 \end{bmatrix}.
\end{equation*}
The determinant of the matrix is 
\begin{equation*}
 \cos(m\theta_1)\sin(m\theta_2)-\cos(m\theta_2)\sin(m\theta_1)
 =\sin(m(\theta_2-\theta_1))\ne 0
\end{equation*}
by $\theta_1-\theta_2\notin \pi\mathbb{Q}$ and $m\ne 0$. 
Hence we have 
\begin{equation*}
 p_{n(j)}=p_{n(j)+1}=0,\ j\in\mathbb{N}^+.
\end{equation*} 
For the case of $m(n(j))=0$, we have  
 \begin{equation*}
 \sum_{\lambda_n=\lambda_j} a_n(z_\ell)p_n
 = \pi^{-1/2}\lambda_j^{-1/2}p_{n(j)}= 0,
\end{equation*}
which gives $p_{n(j)}=0$. Now we have proved $p_n=0$ for $\N+$ and 
the proof is complete.
\end{proof}

\begin{lemma}\label{lemma_uniqueness_Dirichlet}
 Let $\{\tau_n:\N+\}$ be an absolutely convergent complex sequence and 
 $\{\gamma_n:\N+\}$ be a real sequence satisfying 
 $0\le \gamma_1<\gamma_2<\cdots,\ \gamma_n\to \infty.$
For the complex series $\sum_{n=1}^\infty \tau_n e^{-\gamma_n t}$ which 
is defined on $\mathbb{C}^+$, if the set of its zeros on $\mathbb{C}^+$ 
has an accumulation point, then $\tau_n=0,\ \N+$.
\end{lemma}
\begin{proof}
 This lemma can be seen from the analyticity and unique expansion of 
 the generalized Dirichlet series. Here we provide another proof that
makes clear the need for the pieces we have assembled. 
 
 Following the proof of Lemma \ref{analytic}, the analyticity 
 of $e^{-\gamma_n t}$ on $\mathbb{C}^+$ and the absolute 
 convergence of $\{\tau_n:\N+\}$ ensure that 
 $\sum_{n=1}^\infty \tau_n e^{-\gamma_n t}$ is analytic on $\mathbb{C}^+$. 
 Then by the identity theorem for holomorphic functions, if 
 the set of its zeros on $\mathbb{C}^+$ has an accumulation point,
 then $\sum_{n=1}^\infty \tau_n e^{-\gamma_n t}\equiv 0,\ t\in\mathbb{C}^+$. 
 Now we restrict $t$ on $[0,\infty)$ and take Laplace transform.
 By the dominated convergence theorem and the absolute 
 convergence of $\{\tau_n:\N+\}$, we have 
 \begin{equation*}
  0=\L\Big(\sum_{n=1}^\infty \tau_n e^{-\gamma_n t}\Big)(s) 
  =\sum_{n=1}^\infty \tau_n (s+\gamma_n)^{-1},\quad \re s>0.
 \end{equation*}
 From the proof of Lemma \ref{analytic} we can extend the series 
 $\sum_{n=1}^\infty \tau_n (s+\gamma_n)^{-1}$ analytically to 
 $\mathbb{C}\setminus\{-\gamma_n:\N+\}$, consequently, 
 \begin{equation*}
 \sum_{n=1}^\infty \tau_n (s+\gamma_n)^{-1}\equiv 0, 
 \quad s\in\mathbb{C}\setminus\{-\gamma_n:\N+\}.
 \end{equation*}
 Since $\{\gamma_n\}$ is 
 strictly increasing and tends to infinity, it does not contain 
 accumulation points. This means for each $l\in \mathbb{N}^+$, we can 
 take a closed contour which only contains $-\gamma_l$, not 
 $-\gamma_n,\ n\ne l$. Taking integral on both sides of the above 
 equality along this contour, the residue theorem gives that $\tau_l=0$. 
 The proof is complete. 
 \end{proof}

\begin{lemma}\label{lemma_uniqueness_2}
 Given $\epsilon>0$ and the condition $\theta_1-\theta_2\notin \pi\mathbb{Q}$, 
 then 
 \begin{equation*}
  \lim_{\re s\to\infty}e^{\epsilon s} \Big[\sum_{n=1}^\infty a_n(z_\ell)p_n\lambda_n
  (s+\lambda_n)^{-1}\Big]=0,\quad \ell=1,2
 \end{equation*} 
 gives $p_n=0,\ \N+$.
\end{lemma}
\begin{proof}
Fix $\ell\in\{1,2\}$ and define 
\begin{equation*}
F_\ell(t):=
 \begin{cases}
  \sum_{n=1}^\infty a_n(z_\ell)p_n(1-e^{-\lambda_nt}),& t\ge 0;\\
  0,&t<0.
 \end{cases}
\end{equation*}
The Convolution Theorem and \eqref{laplace_analysis} give that 
for $\re s>0$, 
\begin{equation*}
\begin{aligned}
\int_{-\infty}^\infty\! e^{-ts} \int_{-\infty}^\infty \!H(t-\tau+\epsilon) 
F_\ell(\tau)\,d\tau \,dt 
&=\int_{-\infty}^\infty e^{-ts}H(t+\epsilon)\,dt\ \  
\int_{-\infty}^\infty e^{-\tau s} F_\ell(\tau)\,d\tau\\
&=\int_{-\epsilon}^\infty\!e^{-ts}\,dt\ \ 
\int_0^\infty\!e^{-\tau s} \Big[\sum_{n=1}^\infty a_n(z_\ell)p_n
(1-e^{-\lambda_n\tau})\Big]\,d\tau\\
&= s^{-2}e^{\epsilon s} \Big[\sum_{n=1}^\infty a_n(z_\ell)p_n\lambda_n 
(s+\lambda_n)^{-1}\Big],
\end{aligned}
\end{equation*}
which together with the assumption implies that 
\begin{equation*}
 \lim_{\re s\to \infty}\int_{-\infty}^\infty e^{-ts} \int_{-\infty}^\infty H(t-\tau+\epsilon) 
F_\ell(\tau)\,d\tau \,dt=0.
\end{equation*}
A direct calculation then gives 
\begin{equation*}
\begin{aligned}
{\mathbb L}_1 :=& \int_{-\infty}^\infty\! e^{-ts}\! \int_{-\infty}^\infty\! H(t\!-\!\tau+\epsilon) 
F_\ell(\tau)\,d\tau \,dt
=\int_{-\epsilon}^\infty\! e^{-ts} \int_0^{t+\epsilon} 
\sum_{n=1}^\infty a_n(z_\ell)p_n(1-e^{-\lambda_n\tau})\,d\tau \,dt\\
=&\int_{-\epsilon}^\infty\! e^{-ts}\ \Big[\sum_{n=1}^\infty a_n(z_\ell)p_n
(t+\epsilon-\lambda_n^{-1}+\lambda_n^{-1}e^{-\lambda_n(t+\epsilon)})\Big]\,dt\\
=&\int_{-\epsilon}^0\! e^{-ts}\ \Big[\sum_{n=1}^\infty a_n(z_\ell)p_n
(t+\epsilon-\lambda_n^{-1}+\lambda_n^{-1}e^{-\lambda_n(t+\epsilon)})\Big]\,dt\\
&+\int_0^\infty\! e^{-ts}\ \Big[\sum_{n=1}^\infty a_n(z_\ell)p_n
(t+\epsilon-\lambda_n^{-1}+\lambda_n^{-1}e^{-\lambda_n(t+\epsilon)})\Big]\,dt\\
=&\int_0^\epsilon\! e^{(\epsilon-t)s}\ \Big[\sum_{n=1}^\infty a_n(z_\ell)p_n
(t-\lambda_n^{-1}+\lambda_n^{-1}e^{-\lambda_nt})\Big]\,dt\\
&+\int_0^\infty\! e^{-ts}\ \Big[\sum_{n=1}^\infty a_n(z_\ell)p_n
(t+\epsilon-\lambda_n^{-1}+\lambda_n^{-1}e^{-\lambda_n(t+\epsilon)})\Big]\,dt\\
:=&S^\ell_1(s)+S^\ell_2(s),
\end{aligned}
\end{equation*}
where the second equality comes from the absolute convergence of 
$\sum_{n=1}^\infty a_n(z)p_n$ stated by Lemma \ref{convergence}, 
and the term by term calculation. For $S^\ell_2(s)$, 
with the absolute convergence of $\sum_{n=1}^\infty a_n(z)p_n$ and 
\eqref{laplace_analysis}, the summation and integral can be exchanged 
and this leads to the following asymptotic result  
\begin{equation*}
  S^\ell_2(s)=\sum_{n=1}^\infty a_n(z_\ell)p_n \Big[s^{-2}+\epsilon s^{-1}-\lambda_n^{-1}s^{-1}
+\lambda_n^{-1}e^{-\lambda_n\epsilon}(s+\lambda_n)^{-1}\Big]\to 0, 
\quad \re s\to \infty.
\end{equation*}
Now we have 
\begin{equation*}
 \lim_{\re s\to \infty} S^\ell_1(s) =\lim_{\re s\to \infty} 
 {\mathbb L}_1-\lim_{\re s\to \infty} S^\ell_2(s)=0.
\end{equation*}
This implies that $S^\ell_1(s)$ is bounded on $\mathbb{C}^+.$ 
For $s$ with $\re s< 0$, using the fact that
$0<\lambda_1\le\cdots\le\lambda_n\le\cdots\to \infty$, we have 
\begin{equation*}
\begin{aligned}
 |S^\ell_1(s)|&\le \int_0^\epsilon |e^{s(\epsilon-t)}|
 \ \Big[\sum_{n=1}^\infty |a_n(z_\ell)p_n|\ 
|t-\lambda_n^{-1}+\lambda_n^{-1}e^{-\lambda_nt}|\Big]\,dt\\
&\le \int_0^\epsilon C\sum_{n=1}^\infty |a_n(z_\ell)p_n|\,dt<\infty. 
\end{aligned}
\end{equation*}
Hence, we are able to extend the domain of $S^\ell_1(s)$ to the whole 
complex plane $\mathbb{C}$ and its boundedness can be derived.  
By the Cauchy-Riemann equations, it is not hard to show that 
$S^\ell_1(s)$ is holomorphic on $\mathbb{C}$. Namely, $S^\ell_1(s)$ is an 
entire function. The boundedness and Liouville's theorem yield that 
$S^\ell_1\equiv C$ on $\mathbb{C}$, and the limit result means that 
$S^\ell_1\equiv 0$ on $\mathbb{C}$.
Now we have 
\begin{equation*}
\begin{aligned}
 \int_0^\epsilon e^{s(\epsilon-t)}\ \Big[\sum_{n=1}^\infty a_n(z_\ell)p_n&
(t-\lambda_n^{-1}+\lambda_n^{-1}e^{-\lambda_nt})\Big]\,dt\\
&=e^{s\epsilon}\int_0^\epsilon e^{-st}\ \Big[\sum_{n=1}^\infty a_n(z_\ell)p_n
(t-\lambda_n^{-1}+\lambda_n^{-1}e^{-\lambda_nt})\Big]\,dt\equiv 0,
\end{aligned}
\end{equation*}
which means for $\re s>0$, 
\begin{equation*}
\begin{aligned}
 0\equiv&\int_0^\epsilon e^{-st}\ \Big[\sum_{n=1}^\infty a_n(z_\ell)p_n
(t-\lambda_n^{-1}+\lambda_n^{-1}e^{-\lambda_nt})\Big]\,dt\\
= &\int_0^\infty e^{-st}H(\epsilon-t)\ \Big[\sum_{n=1}^\infty a_n(z_\ell)p_n
(t-\lambda_n^{-1}+\lambda_n^{-1}e^{-\lambda_nt})\Big]\,dt \\
=&\L \Big(H(\epsilon-t)\ \big[\sum_{n=1}^\infty a_n(z_\ell)p_n
(t-\lambda_n^{-1}+\lambda_n^{-1}e^{-\lambda_nt})\big]\Big)(s).
\end{aligned}
\end{equation*}
It follows that  
\begin{equation*}
 \sum_{n=1}^\infty a_n(z_\ell)p_n(t-\lambda_n^{-1}+\lambda_n^{-1}e^{-\lambda_nt})
 = 0,\quad t\in(0,\epsilon).
\end{equation*}
By Lemma \ref{convergence}, we can calculate the derivative of 
the above series by termwise differentiation, which gives  
\begin{equation*}
 \sum_{n=1}^\infty a_n(z_\ell)p_n(1-e^{-\lambda_nt}) 
 = \sum_{n=1}^\infty a_n(z_\ell)p_n -\sum_{n=1}^\infty a_n(z_\ell)p_ne^{-\lambda_nt}
 =0, \quad t\in(0,\epsilon).
\end{equation*}
We can see for the above series, the conditions of Lemma 
\ref{lemma_uniqueness_Dirichlet} are satisfied. 
Hence, recalling that $\{\lambda_j:j\in\mathbb{N}^+\}$ is the set of 
distinct eigenvalues, we have 
\begin{equation*}
 \sum_{\lambda_n=\lambda_j}a_n(z_\ell)p_n=0,\quad j\in\mathbb{N}^+,  
 \ \ell=1,2.
\end{equation*}
Now Lemma \ref{lemma_uniqueness_1} allows us to deduce that $p_n=0,\ \N+$  
and completes the proof.
\end{proof}

\subsection{Proof of Theorem \ref{thm:uniqueness}}
Now we are in the position to show the main theorem,  
Theorem \ref{thm:uniqueness}.  
\begin{proof}[Proof of Theorem \ref{thm:uniqueness}]
 Write $q$ and $\tilde{q}$ as  
 \begin{equation*}
  q(t)=\sum_{k=1}^{K} q_kH(t-c_k),\ \tilde{q}(t)=\sum_{k=1}^{\tilde{K}} 
  \tilde{q}_kH(t-\tilde{c}_k), \ t\in[0,\infty),
 \end{equation*}
 and define 
 \begin{equation*}
 \begin{aligned}
  & p_n=\l p(\cdot),\varphi_n(\cdot)\ro,
 && P_\ell(s)= \sum_{n=1}^\infty a_n(z_\ell)p_n\lambda_n(s+\lambda_n)^{-1},\\
  &\tilde{p}_n=\l \tilde{p}(\cdot),\varphi_n(\cdot)\ro,
  &&\tilde{P}_\ell(s)= \sum_{n=1}^\infty a_n(z_\ell)\tilde{p}_n\lambda_n(s+\lambda_n)^{-1},
  \ \ell=1,2.
  \end{aligned}
 \end{equation*}
 Also, denote the infimum of the mesh size of $\{c_k\}$ and 
 $\{\tilde{c_k}\}$ as $\eta$ and $\tilde{\eta}$, respectively. 
With \eqref{laplace}, Lemma \ref{analytic} and the analytic continuation, 
it follows that   
 \begin{equation}\label{equality_2}
  \Big[\sum_{k=1}^Kq_ke^{-c_ks}\Big] P_\ell(s)
  =\Big[\sum_{k=1}^{\tilde{K}}\tilde{q}_ke^{-\tilde{c}_ks}\Big] 
  \tilde{P}_\ell(s),
  \quad s\in\mathbb{C}^+,\ \ell=1,2.
 \end{equation}
 
Now we prove $c_1=\tilde{c}_1$ by contradiction. Assume not, without 
loss of generality, we can set $c_1<\tilde{c}_1$. Then there exists 
$\epsilon>0$ such that 
$\epsilon<\min\{\tilde{c}_1-c_1,\eta\}$, and by multiplying 
$e^{(c_1+\epsilon)s}$ on both sides of \eqref{equality_2} 
we obtain that for $s\in\mathbb{C}^+,\ \ell=1,2$, 
\begin{equation}\label{equality_3}
e^{\epsilon s} q_1 P_\ell(s)=
 -\Big[\sum_{k=2}^Kq_ke^{(c_1-c_k+\epsilon)s}\Big] P_\ell(s)
 +\Big[\sum_{k=1}^{\tilde{K}}\tilde{q}_ke^{(c_1-\tilde{c}_k+\epsilon)s}
 \Big] \tilde{P}_\ell(s).
\end{equation}
The assumption $q\in L^1(0,\infty)$ gives that $K\ge 2$, so that the 
first series in the right side is well defined.  
Since $\re s\ge0$, we have 
\begin{equation*}
 |P_\ell(s)|\le \sum_{n=1}^\infty |a_n(z_\ell)p_n|<\infty,\quad 
 |\tilde{P}_\ell(s)|\le \sum_{n=1}^\infty |a_n(z_\ell)\tilde{p}_n|<\infty, 
\end{equation*}
and considering Remark \ref{condition_q}, it follows that 
\begin{equation*}
 \begin{aligned}
  |q_ke^{(c_1-c_k+\epsilon)s}|&\le C\|q\|_{L^1(0,\infty)} 
  e^{[(1-k)\eta+\epsilon]\re s},\\
  |\tilde{q}_ke^{(c_1-\tilde{c}_k+\epsilon)s}|&\le C\|\tilde{q}\|_{L^1(0,\infty)} 
  e^{[(1-k)\tilde{\eta}+c_1-\tilde{c}_1+\epsilon]\re s}.
 \end{aligned}
\end{equation*}
From the result $\epsilon<\min\{\tilde{c}_1-c_1,\eta\}$ we have 
$-\eta+\epsilon<0,\ c_1-\tilde{c}_1+\epsilon<0$.
These properties give that 
\begin{equation*}
\begin{aligned}
\lim_{\re s\to \infty} \sum_{k=2}^Ke^{[(1-k)\eta+\epsilon]\re s} 
&=\lim_{\re s\to \infty}\frac{1-e^{-(K-1)\eta \re s}}
{1-e^{-\eta \re s}}e^{(-\eta+\epsilon)\re s}=0,\\
\lim_{\re s\to \infty}\sum_{k=1}^{\tilde{K}} 
e^{[(1-k)\tilde{\eta}+c_1-\tilde{c}_1+\epsilon]\re s}&=
\lim_{\re s\to \infty}\frac{1-e^{-\tilde{K}\tilde{\eta} \re s}}
{1-e^{-\tilde{\eta} \re s}}e^{(c_1-\tilde{c}_1+\epsilon)\re s}=0.
\end{aligned}
\end{equation*}
Hence, the right side of \eqref{equality_3} converges to $0$ as 
$\re s\to\infty$, so does the left side, namely  
\begin{equation*}
 \lim_{\re s\to \infty}e^{\epsilon s} q_1 P_\ell(s)=0, \ \ell=1,2.
\end{equation*}
With Lemma~\ref{lemma_uniqueness_2} and the fact $q_1\ne 0$ from 
Assumption~\ref{assumption}, we have 
$p_n=0,\ \N+$. This means $p=0$ in $L^2(\Omega)$ and contradicts 
with Assumption \ref{assumption}. Hence, we have $c_1=\tilde{c}_1$.

Inserting this into \eqref{equality_3} and 
the following equality can be derived
\begin{equation*}
 e^{\epsilon s} \big[q_1 P_\ell(s)-\tilde{q}_1 \tilde{P}_\ell(s)\big]
 =-\Big[\sum_{k=2}^Kq_ke^{(c_1-c_k+\epsilon)s}\Big] P_\ell(s)
 +\Big[\sum_{k=2}^{\tilde{K}}\tilde{q}_ke^{(c_1-\tilde{c}_k+\epsilon)s}\Big] 
 \tilde{P}_\ell(s).
\end{equation*}
Setting $0<\epsilon<\min\{\eta,\tilde{\eta}\}$ and using the above 
limit analysis give that 
the left side of the above equality tends to $0$ as $\re s\to \infty$. 
Now Lemma~\ref{lemma_uniqueness_2} shows that 
$q_1p_n-\tilde{q}_1\tilde{p}_n=0$ for $\N+$.
This means that 
$$\l q_1p(\cdot)-\tilde{q}_1\tilde{p}(\cdot),\varphi_n(\cdot)\ro=0,\quad \N+,$$ 
which together with the completeness of $\{\varphi_n:\N+\}$ in $L^2(\Omega)$
gives $q_1p(x)=\tilde{q}_1\tilde{p}(x)$ in $L^2(\Omega)$. Since 
$q_1, \tilde{q}_1$ are not zero, we can define $C_0:=\tilde{q}_1/q_1$ 
and obviously $C_0\ne 0$. Then we have 
$C_0\,q_1=\tilde{q}_1$ and $p=C_0\tilde{p}$ in $L^2(\Omega)$.

The result $p=C_0\,\tilde{p}$ in $L^2(\Omega)$ implies that 
$P_\ell(s)=C_0\,\tilde{P}_\ell(s)$.
Now, we want to show $C_0\,q(t)=\tilde{q}(t)$. 
Subtracting $q_1e^{-c_1s}P_\ell(s)$ from both sides of \eqref{equality_2} 
gives that 
\begin{equation*}
 \Big[\sum_{k=2}^{K}q_ke^{-c_ks}\Big] P_\ell(s)
  =\Big[\sum_{k=2}^{\tilde{K}}\tilde{q}_ke^{-\tilde{c}_ks}\Big] \tilde{P}_\ell(s).
\end{equation*}
Using the above argument we can obtain $c_2=\tilde{c}_2$ and 
$C_0\,q_2=\tilde{q}_2$. If $K, \tilde{K}$ are both infinity, we can continue 
this procedure and obtain   
\begin{equation*}
 c_k=\tilde{c}_k,\ C_0q_k=\tilde{q}_k,\ k\in \mathbb{N}^+,
\end{equation*}
which means $C_0\,q=\tilde{q}$ on $[0,\infty)$. If the claim that $K=\infty$ 
and $\tilde{K}=\infty$ is not valid, without loss of generality, we 
can assume $K<\infty$. For the case of $K<\tilde{K}$, 
following the above procedure we can get 
\begin{equation}\label{equality_4}
 c_k=\tilde{c}_k,\ C_0\,q_k=\tilde{q}_k,\ k=1,\cdots,K.
\end{equation}
Subtracting $\Big[\sum_{k=1}^{K}q_ke^{-c_ks}\Big] P_\ell(s)$ from 
both sides of \eqref{equality_2}, the following equality can be deduced 
\begin{equation*}
 \Big[\sum_{k=K+1}^{\tilde{K}}\tilde{q}_k\,e^{-\tilde{c}_ks}\Big] 
 \tilde{P}_\ell(s)=0, \quad s\in\mathbb{C}^+,\ \ell=1,2.
\end{equation*} 
This result means that the union of the sets of zeros of 
$\sum_{k=K+1}^{\tilde{K}}\tilde{q}_k\,e^{-\tilde{c}_ks}$ and 
$\tilde{P}_\ell(s)$ should cover $\mathbb{C}^+$. 
The proof of Lemma \ref{lemma_uniqueness_Dirichlet} and the 
condition $\tilde{q}_k\ne 0$ give that 
the set of zeros of $\sum_{k=K+1}^{\tilde{K}}\tilde{q}_k\,e^{-\tilde{c}_ks}$ 
on $\mathbb{C}^+$ does not contain accumulation points, 
so we can find an open connected nonempty subset $\mathbb{C}_1\subset \mathbb{C}^+$ such that 
$\tilde{P}_\ell(s)\equiv0$ on $\mathbb{C}_1$, $\ell=1,2$. Then the analyticity of 
$\tilde{P}_\ell(s)$  
supported by Lemma \ref{analytic} gives that $\tilde{P}_\ell,\ \ell=1,2$ 
vanish on $\mathbb{C}^+$. This together with Lemma 
\ref{lemma_uniqueness_2} leads to $p=\tilde{p}=0$ in $L^2(\Omega)$, which 
contradicts with Assumption \ref{assumption}. 
Similarly, we can derive an analogous contradiction for 
the case of $K>\tilde{K}$. Now we conclude that $K=\tilde{K},$ which  
together with \eqref{equality_4} implies $C_0\,q(t)=\tilde{q}(t)$. 
The proof is complete.
\end{proof}

\begin{remark}\label{generalization}
While we have set this problem in the unit disc and the underlying
elliptic operator  is the negative Laplacian,
 the above proof of uniqueness goes through for an arbitrary domain
$\Omega\subset\mathbb{R}^2$
with smooth boundary and a self-adjoint elliptic operator 
$\mathbb{L} =- \nabla\cdot(a\nabla u) + q u$ where $a(x)\geq a_0>0$ and $q\geq 0$
and with $a$, $q \in L^\infty(\Omega)$.
The essential observation is that the eigenfunctions $\{\varphi_n\}$
form a complete basis for $L^2(\Omega)$ as does their restrictions to
$\partial\Omega$. The latter claim of completeness follows from 
the uniqueness of the Dirichlet problem on $\Omega$.
In addition, the eigenvalues obey the identical asymptotic behavior
as for the negative Laplacian due to Weyl's formula.
This is crucial for the lemmas of this section.
Of course, the statement of Theorem \ref{thm:uniqueness} must now be modified
so as to choose the boundary measurement points $z_{\ell}$ to not coincide
with a zero of any $\varphi_n(x)$ when $x\in\partial\Omega$.
\end{remark}

\section{Numerical reconstruction}
In this section we show numerical reconstructions of $p$ and $q$ from 
boundary flux data measurements following the algorithm described in the proof
of Theorem~\ref{thm:uniqueness}.
In keeping with a practical situation, 
truncated time-value measurements are taken over a finite interval --
in this case $[0,T]$ is used with $T=1$.
We remark that this is actually a long time period as the traditional scaling
of the parabolic equation to unit coefficients means that the diffusion
coefficient $d$ is absorbed into the time variable and our value of $T$
represents the product of the actual final time of measurement and the value
of $d$.
In fact, $d$ is itself the ratio of the conductivity and specific heat.
Values of $d$ of course vary widely with the material but metals for example
have a range of around $10^{-4}$ to $10^{-5}$meters${^2}$/second.

\subsection{Iterative scheme}
For $(\cos{\theta_\ell},\sin{\theta_\ell})\in \partial\Omega$, from 
Corollary \ref{data_formula} and the convergence result 
Lemma \ref{convergence}, we have the following flux representation 
using termwise differentiation, 
\begin{equation}\label{eqn:sol_rep_recon}
  \frac{\partial u}{\n}(1,\theta_\ell,t)=-\sum_{n=1}^\infty 
  a_n(z_\ell)\lambda_n p_n \int_0^t e^{-\lambda_n(t-\tau)}q(\tau)\ d\tau
\end{equation}
where we have again used polar coordinates.
Since the unknown function $p$ is represented by its Fourier 
coefficients $\{p_n\}$, we consider to reconstruct $(p,q)$ in the 
space $\S_N\times L^2[0,T]$, where 
\begin{equation*}
 \S_N=\s\{\varphi_n(x):n=1,\cdots,N\}.
\end{equation*}

We define the forward operator $F$ as 
\vskip-25pt
\begin{equation*}
\qquad\qquad\qquad\qquad
F(p,q)=
 \begin{bmatrix}
  \frac{\partial u}{\n}(1,\theta_1,t)\\\\
  \frac{\partial u}{\n}(1,\theta_2,t)
 \end{bmatrix}
\end{equation*}
and build an iteration scheme to solve 
\vskip-25pt
\begin{equation*}
\qquad\qquad\qquad\qquad
 F(p,q)= g^\delta(t):=\begin{bmatrix}
  g_1^\delta(t)\\\\
  g_2^\delta(t)
 \end{bmatrix}.
\end{equation*}
Here $g^\delta$ is the perturbed measurement satisfying 
$\|(g^\delta-g)/g\|_{C[0,T]}\le \delta$.  
Clearly, if either of $p(x)$ and $q(t)$ is fixed, the operator 
$F$ is linear.
Consequently, we can construct the sequential iteration scheme 
using Tikhonov regularization as 
\begin{equation}\label{iteration}
\begin{aligned}
p_{j+1}:=&\argmin_{p\in \S_N} \|F[q_j]p-g^\delta\|_{L^2(\Omega)}^2+\beta_p\|p\|_{L^2(\Omega)}^2,\\
q_{j+1}:=&\argmin_{q\in L^2[0,T]}\|F[p_j]q-g^\delta\|_{L^2[0,T]}^2+\beta_q\|\triangledown q\|_{L^1[0,T]}.
\end{aligned}
\end{equation}
In the case of $\{q_j\}$, we choose the total variation regularization 
\cite{MuellerSiltanen:2012} to make sure each $q_j$ saves the 
edge-preserving property to fit the exact 
solution $q(t)$, which is a step function. $\beta_p, \beta_q$ are the 
regularizing parameters. 

\subsection{Regularization strategies}
In equation \eqref{eqn:sol_rep_recon} by necessity any use of
this from a numerical standpoint must truncate to a finite sum.
One might be tempted to use ``as many eigenfunctions as possible''
but there are clearly limits imposed by the data measurement process.
Two of these will be discussed in this section.

We will measure the flux at the points $\theta_\ell$ at a series of time steps.
If these steps are $\delta t$ apart, then the exponential term
$e^{-\lambda_n t}$ with $n=N$, the maximum eigenvalue index used,
is a limiting factor: as a multiplier if $e^{-\lambda_N \delta t}$ is
too small relative to the effects caused by any assumed noise in the data,
then we must either reduce $\delta t$ or decrease $N$.
In short, high frequency information can only be obtained from information
arising from very short time measurements.

We also noted that the selection of measurement points $\{\theta_\ell\}$
should be made to avoid zeros of eigenfunctions on the boundary
as otherwise the information coming from these eigenfunctions is unusable.
From the above paragraph, it is clear that only a relatively small number $N$
of these are usable in any event so that we are in fact far from
restricted in any probabilistic sense from selecting the difference in
measurement points even assuming these are all rational numbers when divided
by $\pi$. We can take $\theta =0$ to be the origin of the system without any
loss of generality so that
$\varphi_n(r,\theta)=\omega_n J_m(\sqrt{\lambda_n}r)\{\cos m\theta,\sin m\theta\}$.
If two points at angles $\theta_1$ and $\theta_2$ are taken then
the difference between them is the critical factor; we need to ensure that
$k(\theta_1-\theta_2) \neq j\pi$ for any integers $j,k$.

\begin{figure}[hp!]
	\center
\begin{subfigure}
  \centering
\includegraphics[trim = .9cm .5cm .5cm 4cm, clip=true,height=3.0cm,width=12cm]
		{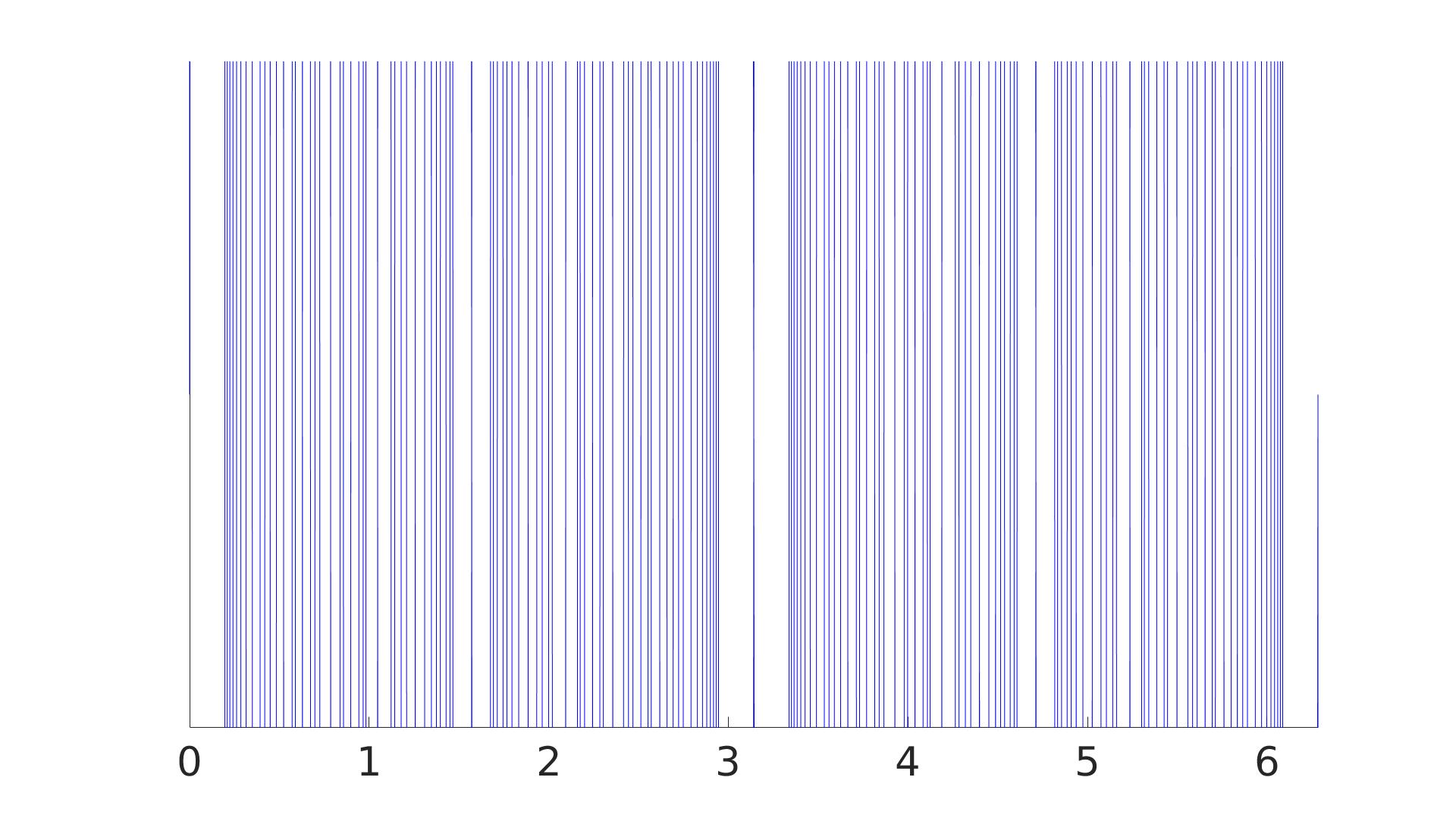}
\end{subfigure}
\caption{\small Gaps between angles.}
\label{gaps}
\end{figure}

Of course the points whose angular difference is a rational number times $\pi$
form a dense set so at face value this might seem a mathematical,
but certainly not a practical, condition.
However, from the above argument,
we cannot use but a relatively small number of eigenfunctions
and so the set of points $(\theta_1,\theta_2)$ with 
$\theta_1-\theta_2\ne(j/k)\pi$
for sufficiently small $k$ might have distinct intervals of sufficient
length for this criteria to be quite practical.
To see this,
consider the rational points generated modulo $\pi$ with denominator
less than the prime value $29$, that is, we are looking for rational
numbers in lowest form $a/b$ with $b<29$ and checking for zeros of
$\sin(a\pi/b)$ for a given $b$.
Clearly taking $b=4$ gives a zero at $\theta=\pi/4$ and we must check
those combinations $a/b$ that would provide a zero close to but less
than $1/4$.
We need only check primes $b$ in the range $2<b<29$ and the fraction closest
to $1/4$ occurs at $a/b=4/17$ which is approximately $0.235$.
Thus the interval that is zero free under this range of $b$ has length
$0.015\pi$ radians or approximately $2.7$ degrees of arc length.
Similar intervals occur at several points throughout the circle.
The gaps in such a situation with $b<29$ is shown in Figure~\ref{gaps}.

Now the question is: if we restrict the eigenvalue index $k$ to be less
than $29$ what range of $m$ index to we obtain and what is the lowest
eigenvalue that exceeds this $k$-range?
Since the $m$-index grows faster than the $k$ for a given eigenvalue index,
we obtain several thousand  eigenvalues,
the largest being approximately $3.5\times 10^4$.
Only with exceedingly small initial time steps we could get such an eigenvalue
and its attendant eigenfunction be utilized in the computations.
If we restrict $k<17$ then the zero-free interval becomes $(\pi/4,4\pi/13)$
with length approximately $10.4$ degrees and
the largest eigenvalue obtained is about $1.5\times 10^4$.
If we decrease down to $k\leq 10$ we get
an angle range of $15.8$ degrees in which to work.

Thus in short, the ill-conditioning of the problem is substantially
due to other factors and not to impossible restrictions on the choice of
observation points $\{\theta_\ell\}$.

\subsection{Numerical experiments}
First we consider the experiment $(e1)$,    
\begin{equation*}
 \begin{aligned}
  (e1):\quad &T=1,
  \ \theta_1=0,\ \theta_2=\frac{13}{32}\pi,\\
  &p(r,\theta)= \frac{5}{\sqrt{30}}\omega_1J_{m(1)}(\sqrt{\lambda_1}r)\cos{(m(1)\theta)}+
  \frac{2}{\sqrt{30}}\omega_2J_{m(2)}(\sqrt{\lambda_2}r)\cos{(m(2)\theta)}\\
  &\qquad\quad\quad+\frac{1}{\sqrt{30}}\omega_2J_{m(2)}(\sqrt{\lambda_2}r)\sin{(m(2)\theta)},\\  
  &q(t)=\chi_{_{[0,1/3)}}+2\chi_{_{[1/3,2/3)}}+1.5\chi_{_{[2/3,1]}}.
  \end{aligned}
\end{equation*}
We use noise-polluted flux measurements on the boundary points at noise
levels ranging from 1\% to 5\% and choose the time measurement step $\delta t$
to be $0.01$.

In order to avoid the loss of accuracy caused by the multiplication between 
$p$ and $q$, we use the normalized exact solution of $p(x)$, namely, 
let $\|p\|_{L^2(\Omega)}=1$. To achieve this setting, in the programming 
of iteration \eqref{iteration}, after each iterative step, we set 
$p_j=p_j/\|p_j\|_{L^2(\Omega)},\ q_j=\|p_j\|_{L^2(\Omega)}q_j$. 
Also, the initial guess $p_0$ and $q_0$ are set as 
\begin{equation*}
\begin{aligned}
 p_0(x)&\equiv 1,\ x\in \Omega,\\
 q_0&:=\argmin_{q\in L^2[0,T]}\|F[p_0]q-g^\delta\|_{L^2[0,T]}^2+\beta_q\|\triangledown q\|_{L^1[0,T]}.
\end{aligned}
\end{equation*}
Depending on the noise level $\delta,$ the values of regularized 
parameters $\beta_p$, $\beta_q$ are 
picked empirically and here the values used are
$\beta_p=1\times10^{-2}$, $\beta_q=8\times10^{-4}$. 
After $j=10$ iterations, the approximations $p_j,q_j$ are recorded and 
displayed by Figure \ref{e1_pq_1}.
This indicates effective numerical convergence of the scheme. 
The errors of approximations 
upon different noise levels are displayed by the following table. 
\begin{center}
\begin{tabular}{|c| c| c|c| }
\hline
                               & $\delta=1\%$ & $\delta=3\%$ &$\delta=5\%$\\
\hline  $\|p-p_j\|_{L^2(\Omega)}$   & $1.34e-1$ & $1.76e-1$& $1.87e-1$\\
\hline  $\|q-q_j\|_{L^2[0,T]}$ & $8.08e-2$ & $8.25e-2$& $9.76e-2$\\
\hline
\end{tabular}
\label{error}
\end{center}
The satisfactory reconstructions shown by the table confirm that 
the iterative scheme \eqref{iteration} is a feasible approach to 
solve this nonlinear inverse problem numerically.  

\setlength{\linewidth}{6.2true in}
\bigskip
\newdimen\xfiglen
\newdimen\yfiglen
\newbox\figurelegend
\newbox\figureone
\font\sevenrm=cmr7
\xfiglen=1.8 true in
\yfiglen=1.2 true in
\setbox\figurelegend=\vbox{\hsize=0.2\xfiglen
\beginpicture
\footnotesize
  \setcoordinatesystem units <0.2\xfiglen,0.1\yfiglen> point at 0 0
  \setplotarea x from 0 to 1, y from 0 to 3
\linethickness=0.8pt
\setsolid
\blue{\relax\putrule from 0.3 2 to 0.6 2 \relax} \put {$q$} [r] at 0.05 2
\setdots <2pt>
\red{\relax\putrule from 0.3 1 to 0.6 1 \relax} \put {$q_j$} [r] at 0.05 1
\endpicture
}
\setbox\figureone=\vbox{\hsize=\xfiglen
\beginpicture
\footnotesize
\sevenrm
  \setcoordinatesystem units <\xfiglen,0.333\yfiglen> 
  \setplotarea x from 0 to 1, y from 0 to 3
  \axis bottom shiftedto y=0 ticks short numbered from 0 to 1 by 0.2 /
  \axis left ticks short numbered from 0 to 3 by 1 /
\linethickness=0.7pt
\putrule from 0 1 to 0.333 1
\setdots <3pt>
\setlinear
\relax
\multiput {\red{\small .}} at 
0    0.98232
0.01 0.98232
0.02 0.98232
0.03 0.98232
0.04 0.98232
0.05 0.98232
0.06 0.98232
0.07 0.98232
0.08 0.98232
0.09 0.98232
0.1  0.98232
0.11 0.98232
0.12 0.98232
0.13 0.98232
0.14 0.98232
0.15 0.98232
0.16 0.98232
0.17 0.98232
0.18 0.98232
0.19 0.98232
0.2  0.98232
0.21 0.98232
0.22 0.98232
0.23 0.98232
0.24 0.98232
0.25 0.98232
0.26 0.98232
0.27 0.98232
0.28 0.98232
0.29 0.98232
0.3  0.98232
0.31 0.98232
0.32 0.98232
0.33 1.5476
0.34 1.8951
0.35 1.8951
0.36 1.8951
0.37 1.8951
0.38 1.8951
0.39 1.8951
0.4  1.8951
0.41 1.8951
0.42 1.8951
0.43 1.8951
0.44 1.8951
0.45 1.8951
0.46 1.8951
0.47 1.8951
0.48 1.8951
0.49 1.8951
0.5  1.8951
0.51 1.8951
0.52 1.8951
0.53 1.8951
0.54 1.8951
0.55 1.8951
0.56 1.8951
0.57 1.8951
0.58 1.8951
0.59 1.8951
0.6  1.8951
0.61 1.8951
0.62 1.8951
0.63 1.8951
0.64 1.8951
0.65 1.8951
0.66 1.8951
0.67 1.8951
0.68 1.5033
0.69 1.5033
0.7  1.5033
0.71 1.5033
0.72 1.5033
0.73 1.5033
0.74 1.5033
0.75 1.5033
0.76 1.5033
0.77 1.5033
0.78 1.5033
0.79 1.5033
0.8  1.5033
0.81 1.5033
0.82 1.4962
0.83 1.4962
0.84 1.4962
0.85 1.4962
0.86 1.4962
0.87 1.4962
0.88 1.4817
0.89 1.4817
0.9  1.4817
0.91 1.4816
0.92 1.4816
0.93 1.4816
0.94 1.4816
0.95 1.4816
0.96 1.4816
0.97 1.4816
0.98 1.4816
0.99 1.4816
1    1.4816 /

\setsolid
\blue{\relax\putrule from 0 1 to 0.333 1 \relax}\relax
\blue{\relax\putrule from 0.334 2 to 0.666 2 \relax}
\blue{\relax\putrule from 0.667 1.5 to 1 1.5 \relax}
\put {\copy\figurelegend} [lt] at 0.05 3
\endpicture
}

\newbox\figuretwo

\setbox\figuretwo=\vbox{\hsize=0.7\linewidth\parindent=0pt
\noindent\includegraphics[trim={7cm 7cm 5cm 
0cm},clip=true,scale=0.17]{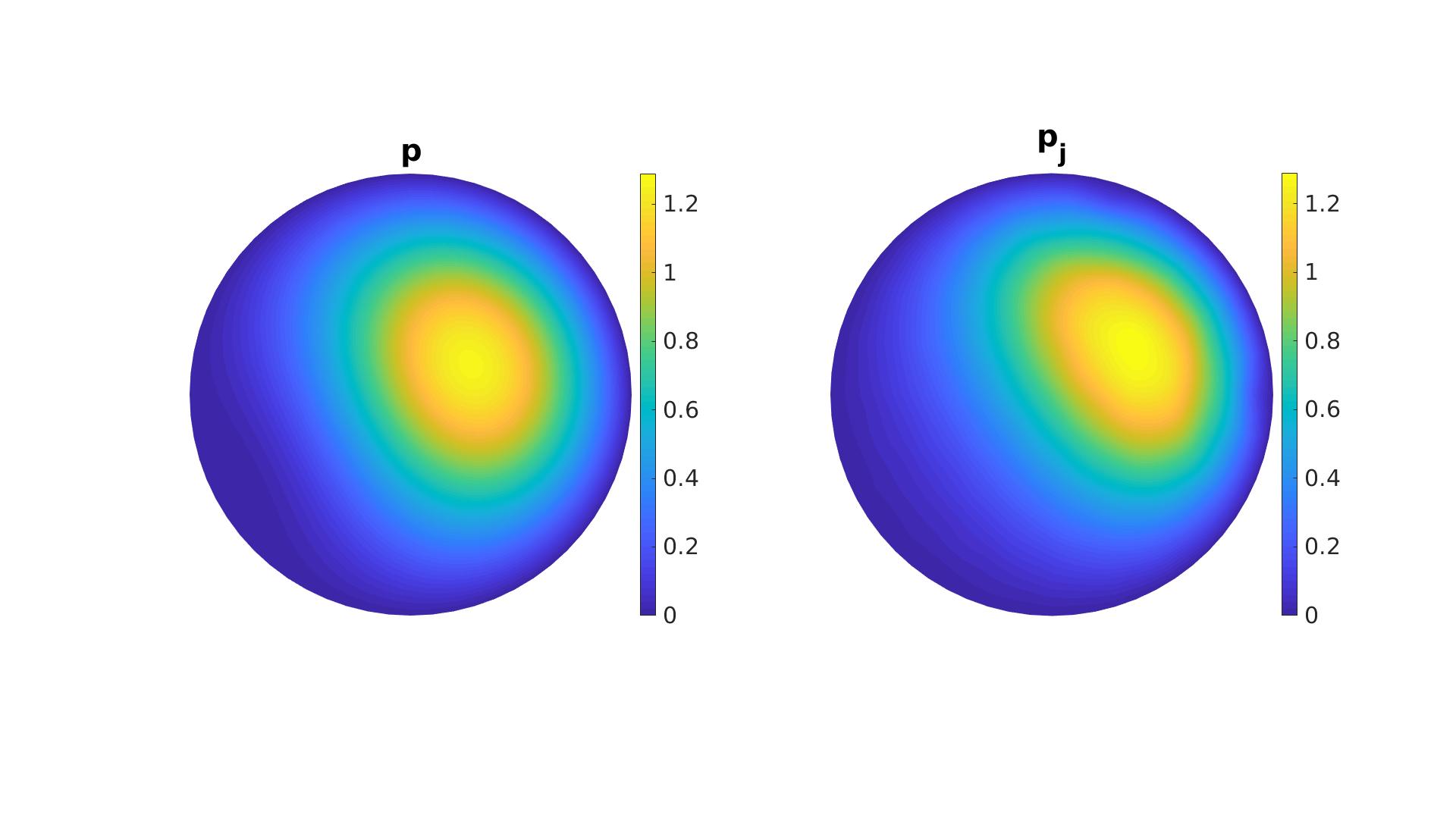}
}
\begin{figure}[ht] 
   \centering
\hbox to\hsize{\copy\figuretwo\hss\qquad\raise3em\copy\figureone}
\caption{\small Experiment $(e1)$, $p$ (left), $p_j$ (center) and 
  $q,\,q_j$ (right). Noise $\delta=1\%$. }
\label{e1_pq_1}
\end{figure}

Next, we seek recovery of a more general $p(x)$:
\begin{equation*}
 \begin{aligned}
  (e2):\quad 
  &p(r,\theta)=\chi_{{}_{r\le 0.5+0.2\cos{2\theta}}} ,\\  
  &q(t)=\chi_{_{[0,1/3)}}+2\chi_{_{[1/3,2/3)}}+1.5\chi_{_{[2/3,1]}},\\
  (e3):\quad 
  &p(r,\theta)=\chi_{{}_{r\le 0.25+0.1\cos{2\theta}}} ,\\  
  &q(t)=\chi_{_{[0,1/3)}}+2\chi_{_{[1/3,2/3)}}+1.5\chi_{_{[2/3,1]}}.
 \end{aligned}
\end{equation*}

\begin{figure}[ht]
	\center
\begin{subfigure}
  \centering
\includegraphics[trim = .5cm 2cm .5cm 2cm, clip=true,height=5.5cm,width=6cm]
		{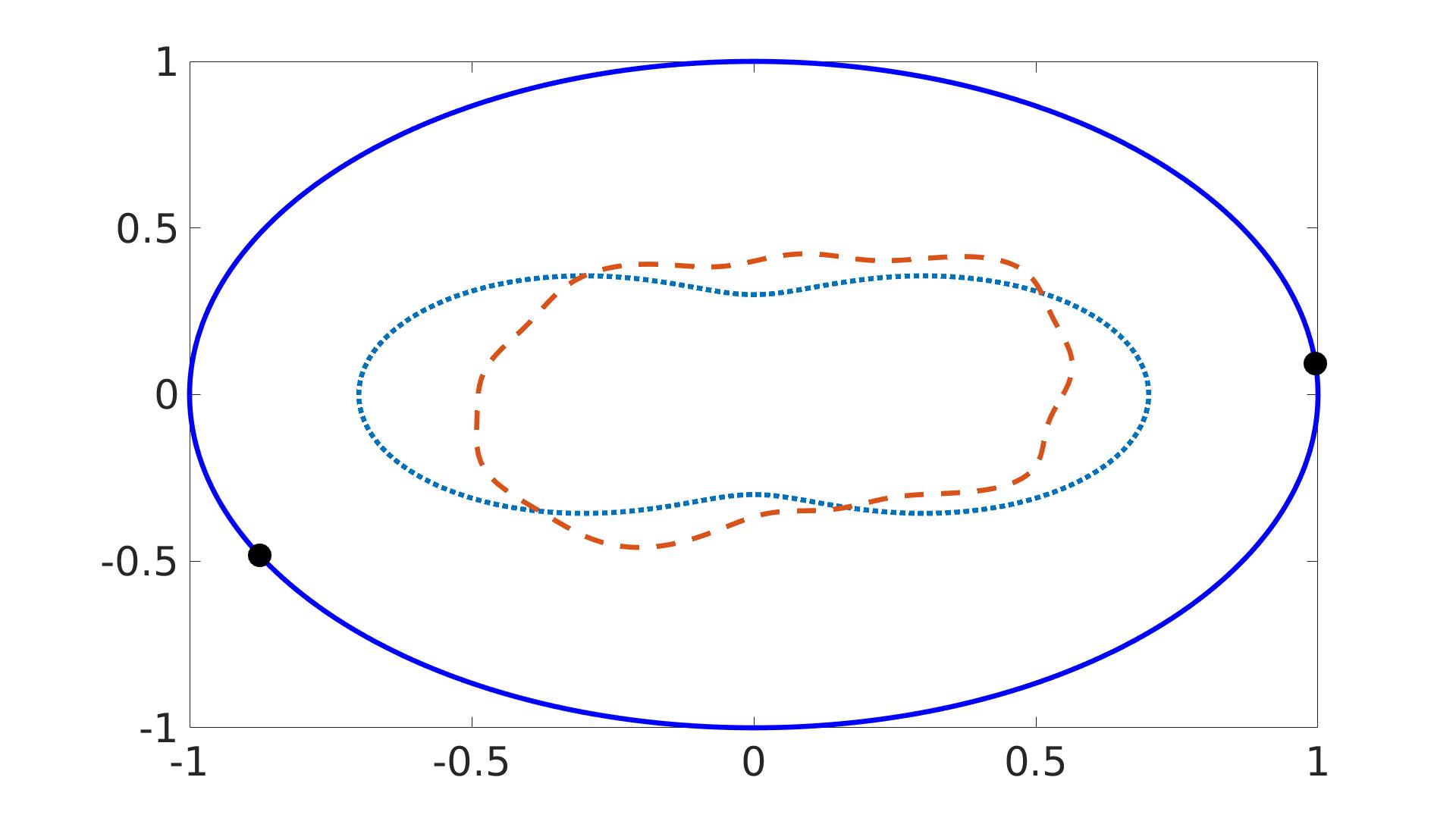}
\includegraphics[trim = .5cm 2cm .5cm 2cm, clip=true,height=5.5cm,width=6.5cm]
		{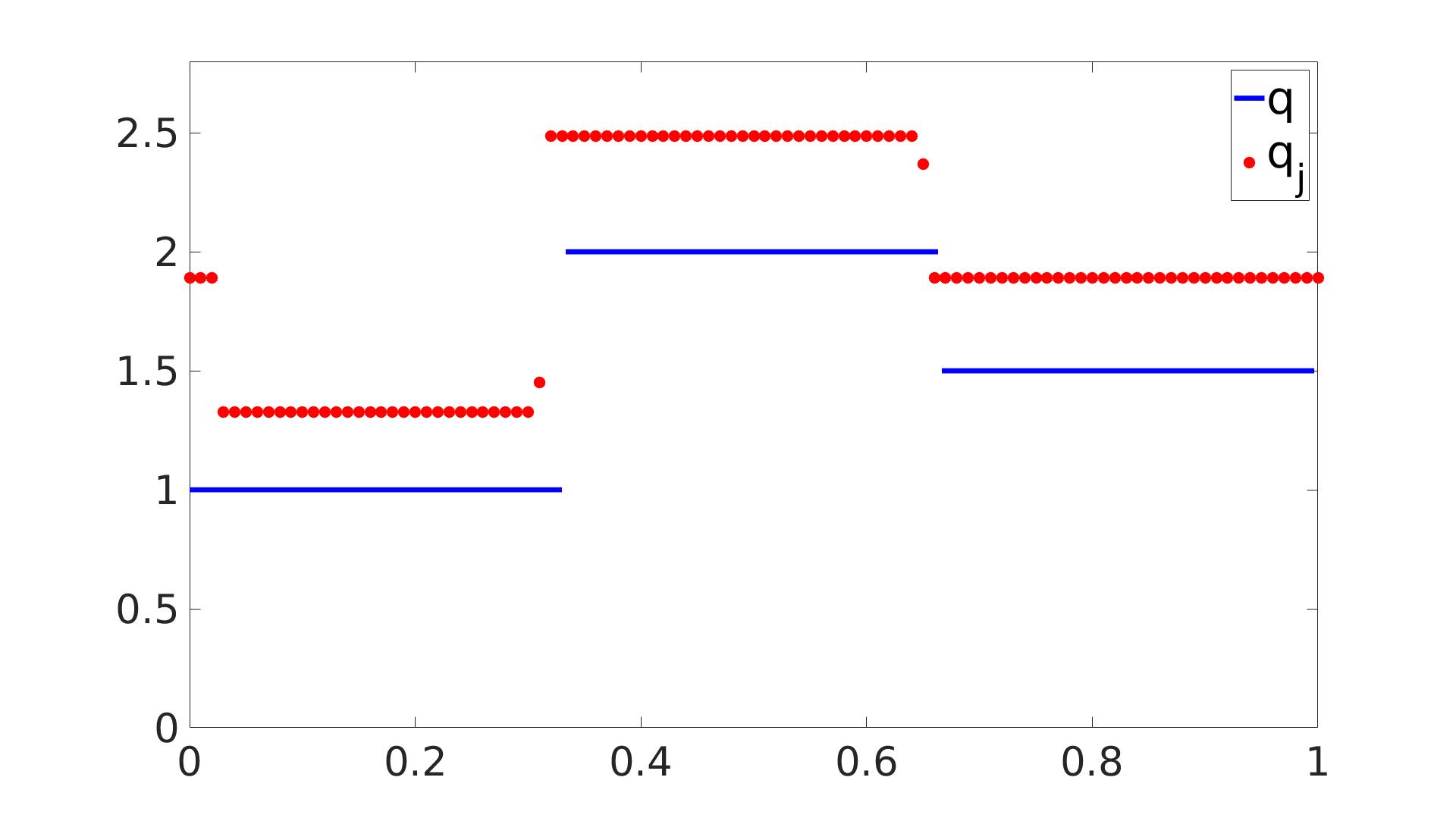}
\end{subfigure}
\caption{\small Experiment $(e2)$, $p$ (left) and $q$ (right), 
$\delta=1\%$.}
\label{e2}
\end{figure}

\begin{figure}[ht]
	\center
\begin{subfigure}
  \centering
\includegraphics[trim = .5cm 2cm .5cm 2cm, clip=true,height=5.5cm,width=6cm]
		{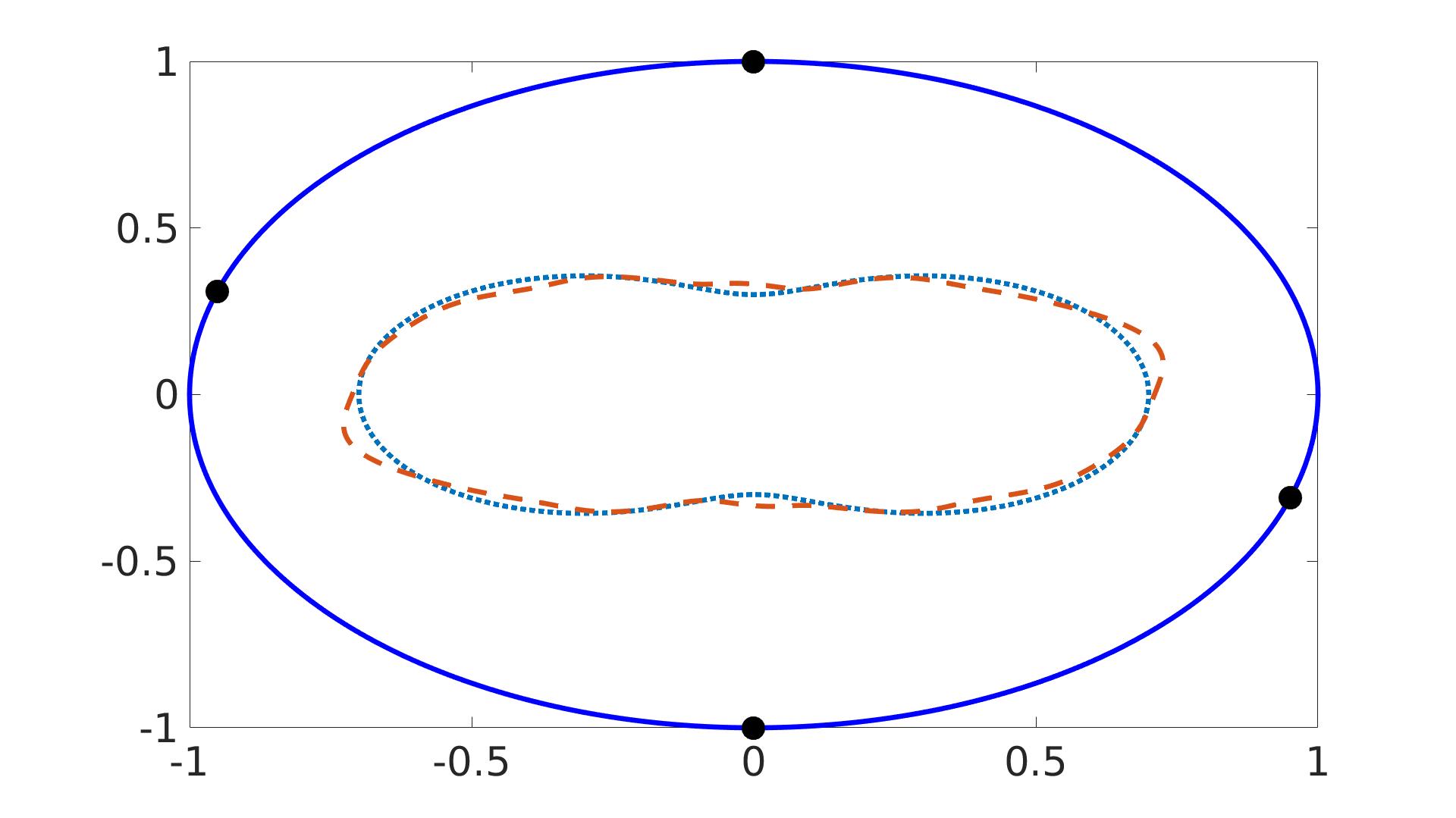}
\includegraphics[trim = .5cm 2cm .5cm 2cm, clip=true,height=5.5cm,width=6.5cm]
		{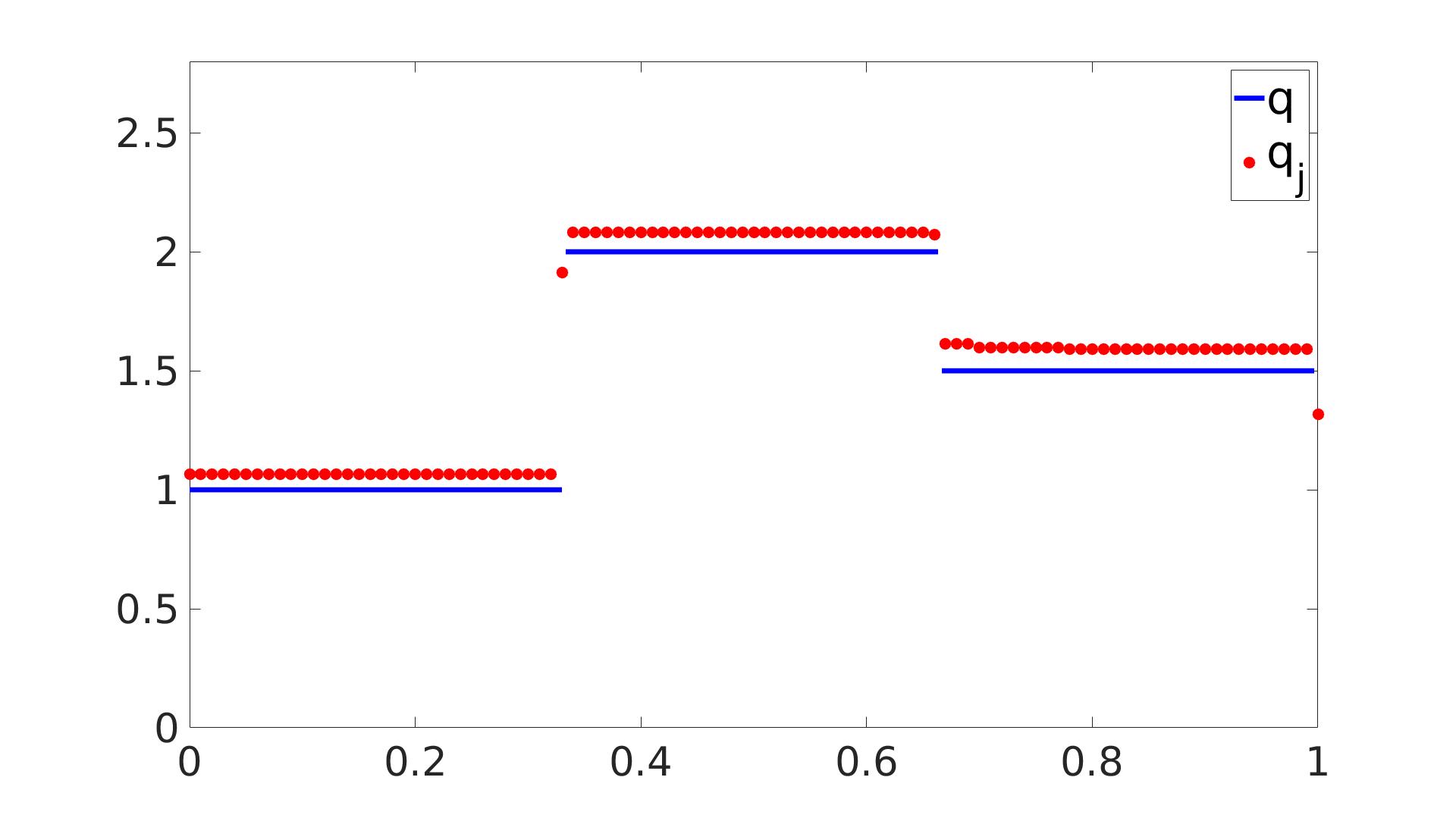}
\end{subfigure}
\caption{\small Experiment $(e2)$ with $4$ measurement points, $p$ (left) and $q$ (right), 
$\delta=1\%$.}
\label{e2_4point}
\end{figure}

\begin{figure}[ht]
	\center
\begin{subfigure}
  \centering
\includegraphics[trim = .5cm 2cm .5cm 2cm, clip=true,height=5.5cm,width=6cm]
		{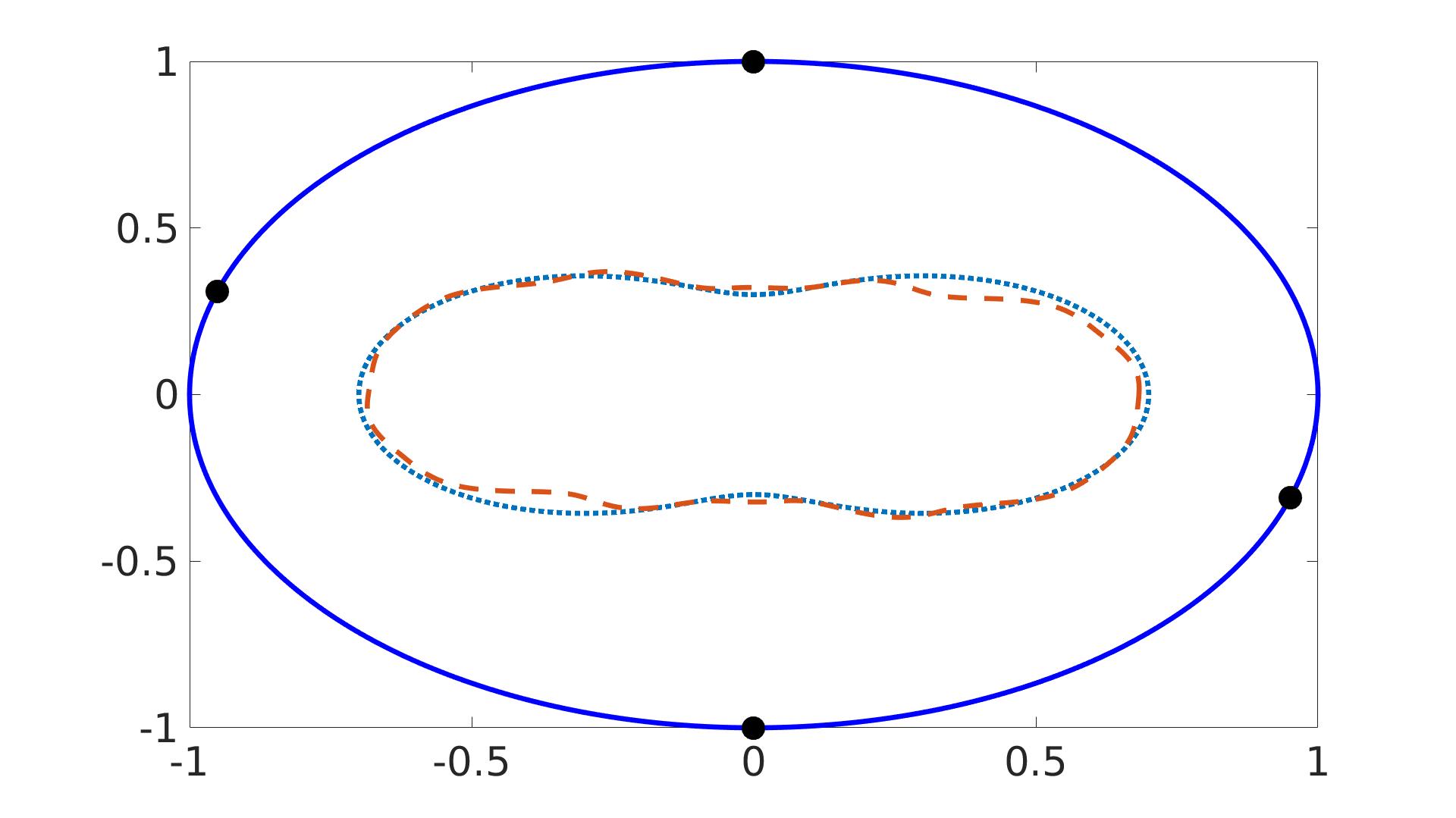}
\includegraphics[trim = .5cm 2cm .5cm 2cm, clip=true,height=5.5cm,width=6.5cm]
		{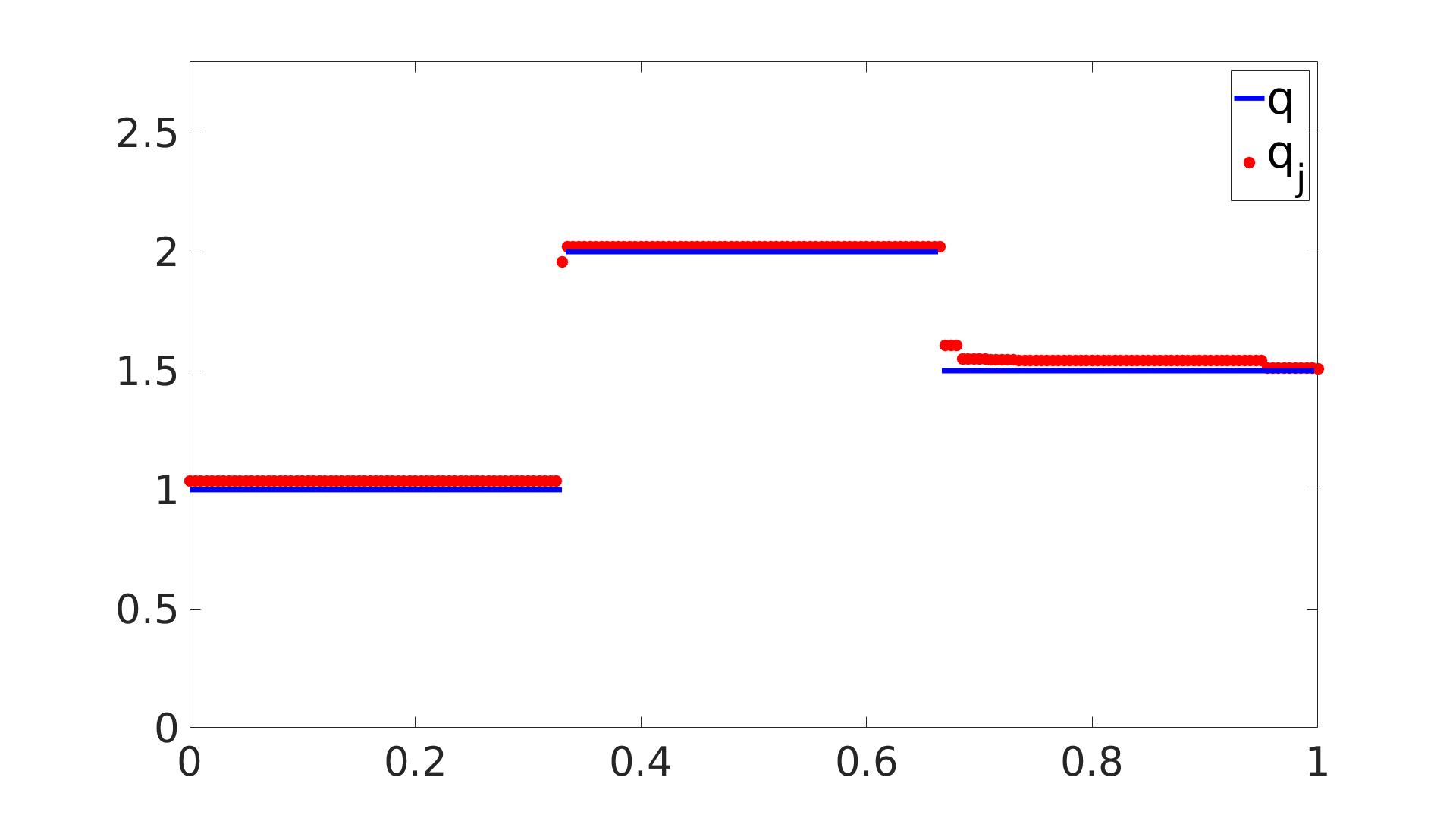}
\end{subfigure}
\caption{\small Experiment $(e2)$ with $4$ measurement points and 
$\delta t=5\times10^{-3}$, $p$ (left) and $q$ (right), 
$\delta=1\%$.}
\label{e2_4point_k200}
\end{figure}

In experiment $(e2)$, a discontinuous, star-like supported exact solution 
$p(x)$ is considered, where the radius function is 
$r(\theta)=0.5+0.2\cos{2\theta}$. We can see this $p$ is out of 
Assumption \ref{assumption}, so the iteration \eqref{iteration} may not 
be appropriate here and in fact, we use the 
Levenberg--Marquardt algorithm to recover the radius function $r(\theta)$, 
see \cite{RundellZhang:2017JCP} for details. 
The numerical results are presented in Figures \ref{e2}, \ref{e2_4point} 
and \ref{e2_4point_k200}, in which the blue dotted line 
and the red dashed line mean the boundaries of $\text{supp}(p)$ and 
$\text{supp}(p_j)$, respectively, and the black bullets are the locations 
of observation points.

Figures~\ref{e2_4point} and \ref{e2_4point_k200} show that with sufficient data, for example, more 
measurement points and finer mesh on time $t$, precise reconstructions 
can be obtained even though Assumption~\ref{assumption} is violated. 
These results indicate that, if we do not pursue the global uniqueness stated by
Theorem~\ref{thm:uniqueness}, 
which requires Assumption~\ref{assumption}, the conditions on 
$p$ and $q$ may be weakened in numerical computations.
This inspires future work on such inverse source problems in order to
provide a rigorous mathematical justification for allowing such inclusions.

\begin{figure}[ht]
	\center
\begin{subfigure}
  \centering
\includegraphics[trim = .5cm 2cm .5cm 2cm, clip=true,height=5.5cm,width=6cm]
		{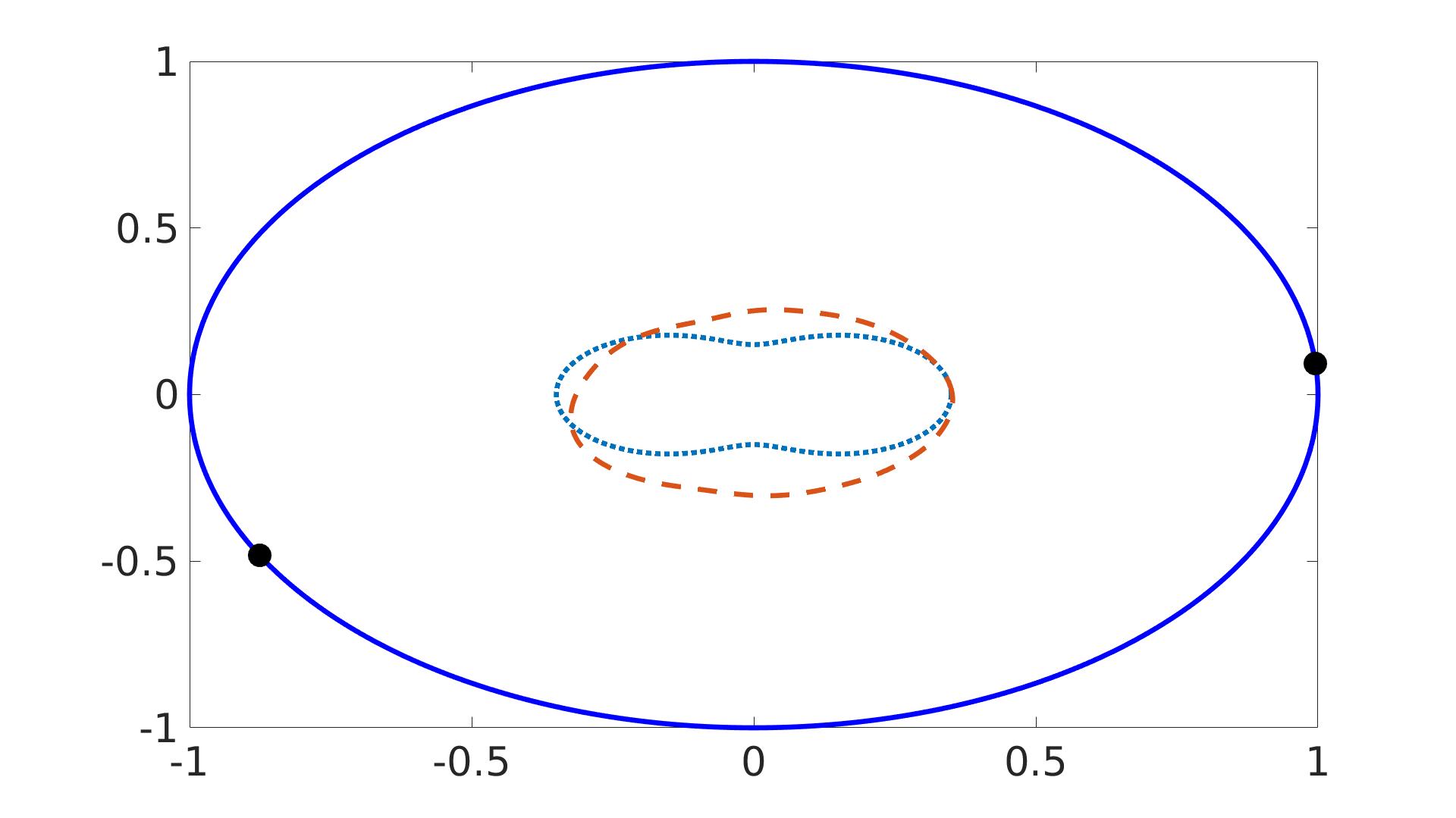}
\includegraphics[trim = .5cm 2cm .5cm 2cm, clip=true,height=5.5cm,width=6.5cm]
		{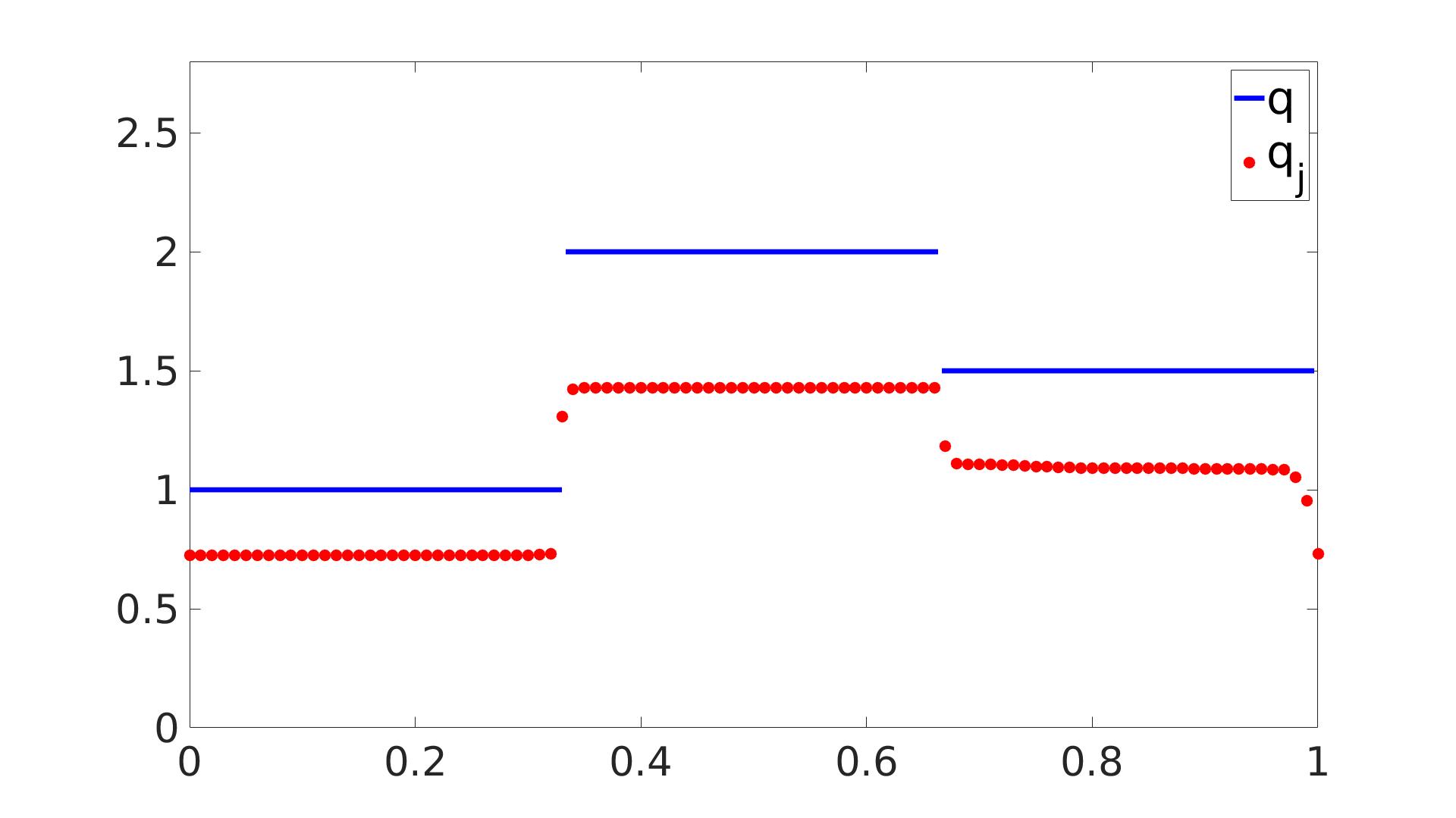}
\end{subfigure}
\caption{\small Experiment $(e3)$, $p$ (left) and $q$ (right), 
$\delta=1\%$.}
\label{e3}
\end{figure}

If we use equation \eqref{eqn:direct_pde} to describe the diffusion of 
pollutants, then $\text{supp}(p)$ means the severely polluted area. With the 
consideration of safety and cost, observations of the flux data 
should be made as far as possible to $\text{supp}(p)$.
This is the reason 
why we set the experiment $(e3)$, in which $p(x)$ has a smaller support. 
Due to the long distance between $\text{supp}(p)$ and the observation 
points, worse results can be expected. See Figure \ref{e3}. 
Hence, accurate and efficient algorithms for this inverse source problem 
with a small $\text{supp}(p)$ are worthy of investigation. 
Of course, in the limit that these become point sources described by
Dirac-delta functions then other tools are available.
See, for example, \cite{HankeRundell:2011}.

\section{Concluding remark and future work}
This paper considers the unique determination of a nonlinear source term 
in the heat equation, which contains two independent unknowns. Only 
finite (here is two) flux measurements are sufficient to support this 
uniqueness, provided some restrictions on $p,q$ stated by Assumption 
\ref{assumption}. Here a natural question may be asked, can we 
weaken the conditions on $p,q$ and meanwhile keep the uniqueness result. 
Let's review the roles of such conditions in the uniqueness proof. 
The smoothness condition $p\in \d$ ensures Lemma \ref{convergence}, 
the absolute convergence of the series $\sum_{n=1}^\infty a_n(z)p_n$, 
which supports the well-definedness of the Laplace transform \eqref{laplace} and 
Lemma \ref{analytic}. While the step function form of $q$ is set for 
the proof of the main theorem, Theorem \ref{thm:uniqueness}. The 
Laplace transform of Heaviside function is the natural exponential function, 
which can not be factored with rational functions, i.e. 
$a_n(z)p_n\lambda_n(s+\lambda_n)^{-1}$. This means we can isolate each 
coefficient pair $(q_k,c_k)$ of $q$ with others in the uniqueness proof 
and then deduce the uniqueness result $p(x)=C_0\tilde{p}(x)$ for the 
space unknown $p$. After this step, the nonlinear inverse problem 
is linearized and naturally, the uniqueness of time unknown $q$ is 
derived. To sum up, to weaken the conditions on $p$ and $q$, 
a new approach may need to be constructed rather than the Laplace 
transform. 

However, the numerical experiments $(e2)$ and $(e3)$ seem to provide a feasible 
way. In the numerical reconstruction aspect, we may consider more 
general unknowns, for instance, discontinuous $p(x)$ and even 
continuous $q(t)$. But in the numerical analysis, we may only prove 
the local uniqueness result, not the global one as Theorem 
\ref{thm:uniqueness}. It may be regarded as the cost for a wider 
class of unknowns. 

Furthermore, extending this work to fractional diffusion equations is 
interesting and meaningful. The fractional case to recover the space 
dependent source $f(x,t)=\chi_{{}_D}$ was considered in  
\cite{RundellZhang:2017JCP}. In the fractional diffusion equation, 
the regular time derivative $\partial/\partial t$ is replaced by the 
fractional derivative $\partial_t^\alpha,\ \alpha\in(0,1).$ 
The fundamental solution for such equations is in terms of Mittag-Leffler 
function $E_{\alpha_1,\alpha_2}(-z)$, but not the natural exponential 
function. This function also holds the analytic property, which means 
the uniqueness proof seems to work. Also, comparing with 
the natural exponential function, the polynomial decay rate of 
$E_{\alpha_1,\alpha_2}(-z)$ may  
cause different performance in the numerical reconstruction. 
In addition, if the fractional order $\alpha$ is set to be unknown, this 
inverse problem will become more challenging.

\section*{Acknowledgment}

\noindent
The work of the first author was supported in part by the
National Science Foundation through award {\sc dms}-1620138.
The second author was supported by Academy of Finland, 
grants 284715, 312110 and the Atmospheric mathematics project of 
University of Helsinki.
\bibliographystyle{abbrv}
\bibliography{rz}

\end{document}